\documentclass[11pt]{article}

\usepackage[margin=1in]{geometry}
\usepackage{amsmath,amssymb,amsthm,mathtools}
\usepackage{enumitem}
\usepackage{booktabs}
\usepackage{array}
\usepackage{float}
\usepackage{xcolor}
\usepackage{tikz}
\usetikzlibrary{arrows.meta}
\usepackage{hyperref}

\hypersetup{
  pdftitle={A Quasisymmetric Analogue of the Bruhat Order},
  pdfauthor={Ryan Bansal},
  colorlinks=true,
  linkcolor=blue!55!black,
  citecolor=blue!55!black,
  urlcolor=blue!55!black
}

\theoremstyle{plain}
\newtheorem{theorem}{Theorem}[section]
\newtheorem{proposition}[theorem]{Proposition}
\newtheorem{lemma}[theorem]{Lemma}
\newtheorem{corollary}[theorem]{Corollary}
\newtheorem{fact}[theorem]{Fact}
\newtheorem{conjecture}[theorem]{Conjecture}

\newtheorem*{thmA}{Theorem A}
\newtheorem*{thmB}{Theorem B}
\newtheorem*{conjC}{Conjecture C}

\theoremstyle{definition}
\newtheorem{definition}[theorem]{Definition}

\newtheorem{observation}[theorem]{Observation}

\theoremstyle{remark}
\newtheorem{remark}[theorem]{Remark}

\newcommand{\ZZ}{\mathbb{Z}}

\newcommand{\supp}{\operatorname{supp}}
\newcommand{\IN}{\operatorname{IN}}
\newcommand{\val}{\operatorname{val}}
\newcommand{\code}{\mathbf{c}}
\newcommand{\rk}{\operatorname{rk}}
\newcommand{\PF}{P_F}

\newcommand{\OmegaF}{\Omega_F}
\newcommand{\rl}{\operatorname{rl}}
\newcommand{\Base}{\operatorname{Base}}
\newcommand{\Spine}{\operatorname{Spine}}
\newcommand{\UD}{\operatorname{UD}}
\newcommand{\cat}{\operatorname{Cat}}
\newcommand{\mex}{\operatorname{mex}}

\newcommand{\forestline}[2]{%
  \draw[gray!65] (#1,0)--(#2,0);
  \foreach \x in {#1,...,#2} {
    \node[gray!70] at (\x,0) {\scriptsize $\times$};
  }
}

\newcommand{\emptyforest}{\varnothing}

\title{A Quasisymmetric Analogue of the Bruhat Order}
\author{Ryan Bansal}
\date{\today}

\begin{document}
\maketitle

\begin{abstract}
  We define a quasisymmetric analogue of the Bruhat order on indexed binary
  forests using the multiplicity-free Monk rule for forest polynomials. For
  the finite window $\mathcal{P}_n$ of forests supported in $[1,n]$, we
  construct an explicit partition into saturated chains terminating at the
  top rank, proving that $\mathcal{P}_n$ is strongly Sperner. We also recover
  the known Catalan ballot formula for its rank numbers.
  
  For principal lower intervals, we prove that $[\varnothing,F]$ is a lattice
  if and only if it is Boolean, and characterize this property by the
  avoidance of a bent corner in $F$, a path of three internal vertices that
  changes direction. This gives a forest analogue of the pattern-avoidance
  characterization of Boolean principal Bruhat intervals. We also conjecture
  that every principal lower interval is shellable, providing computational
  evidence and proving the conjecture for a family of hook intervals.
\end{abstract}

\section{Introduction}
\label{sec:introduction}

Nadeau and Tewari~\cite{NT23} developed forest polynomials in their study of
the cohomology class of the permutahedral variety. They form a basis
of $\ZZ[x_1,x_2,\dots]$ indexed by binary forests~\cite[Proposition~6.7]{NST26}, with nonnegative integral
structure constants~\cite[Theorem~6.3]{NT23}. The forest Monk rule of Nadeau--Spink--Tewari
\cite[Theorem~7.4]{NST26} states that
\[
  (x_1+\cdots+x_i)P_F=\sum_G P_G,
\]
where the summands are distinct. We declare $F\lessdot G$ whenever $P_G$
occurs in such a product, for some $i\geq1$, and call the resulting order the
\emph{quasisymmetric Bruhat order on forests}.

The classical model is the strong Bruhat order on permutations. For the Schubert
polynomials introduced by Lascoux and Sch\"utzenberger~\cite{LS82}, Monk's formula
expresses multiplication by $x_1+\cdots+x_i$ as a sum indexed by Bruhat
covers~\cite[Equation~(4.15$''$)]{Mac91}. The Bruhat order enjoys several
remarkable properties: it is strongly Sperner~\cite[Theorem~3.1]{Sta80}, and
its intervals are Eulerian~\cite{BB05} and shellable~\cite{BW82}.
We ask which of these properties survive for the quasisymmetric Bruhat order.

The algebraic parallel extends beyond multiplication. Forest polynomials
give a basis of the quasisymmetric coinvariant ring~\cite[Theorem~9.7(2)]{NST26},
whose dimension is Catalan rather than factorial~\cite[Theorem~1.1]{ABB04}.
Bergeron, Gagnon, Nadeau, Spink, and Tewari identify this ring with the
cohomology ring of a quasisymmetric flag variety~\cite[Theorem~A]{BGNST26}.
The classical Bruhat order restricted to algebraic noncrossing partitions
also appears in this geometric setting~\cite[Theorem~11.11]{BGNST26}.
This restricted order is isomorphic to the excedance quotient studied by
Bergeron and Gagnon~\cite[Proposition~3.3 and Theorem~4.2]{BG24}.
We introduce a distinct order through Monk multiplication in the forest
basis, prove that its finite windows are strongly Sperner, and
characterize its Boolean principal lower intervals.

There are two natural ways to restrict the quasisymmetric Bruhat order.
The finite window $\mathcal P_n$ consists of forests supported in $[1,n]$;
a principal lower interval $[\emptyforest,F]$ consists of all forests
below one fixed top. The window has $\cat_n$ maximal elements, whereas
the excedance quotient of $S_{n+1}$ has a unique maximum. Thus these
posets are not isomorphic for $n\geq2$, despite both having
$\cat_{n+1}$ elements. Principal lower intervals also exhibit differences
from classical Bruhat intervals: even a three-vertex top can have a
non-Eulerian lower interval.

An \emph{antichain} is a set of pairwise incomparable elements. A ranked
poset is \emph{strongly Sperner} if no union of $r$ antichains is larger than
the union of its $r$ largest rank levels. Our first main result proves this property by an explicit chain partition.

\begin{thmA}[Strong Sperner property; Theorem~\ref{thm:chain-partition} and
Corollary~\ref{cor:strong-sperner}]
The window $\mathcal{P}_n$ is strongly Sperner: it admits an explicit partition
into saturated chains, each ending in the top rank.
\end{thmA}

The construction uses the code $\code(F)=(c_1,c_2,\dots)$ of a forest.
Increasing $c_1$ gives a Monk cover and, below top rank, stays in the window.
Fixing the remaining entries therefore partitions the window into chains,
each ending in rank $n$; every antichain meets each chain at most once.
Figure~\ref{fig:pumping-highways} shows the construction in $\mathcal P_3$.
It also yields sharp bounds for unions of antichains. The ballot rank numbers
give $|\mathcal P_n|=\cat_{n+1}$ and largest antichains of size $\cat_n$.

Chain decompositions also occur in the classical setting: the strong Bruhat
order on permutations admits a symmetric chain decomposition
(an observation Stanley credits to E.~Gansner~\cite[Section~7, p.~182]{Sta80}).
However, for the weak order on permutations, the Sperner
question goes back to Bj\"orner~\cite{Bjo84} in 1984. Stanley's 2017 determinant
conjecture~\cite[Conjecture~2.2]{Sta17} would imply the strong Sperner property.
Gaetz and Gao established that property using a combinatorially defined
$\mathfrak{sl}_2$-action in a 2018 preprint, published in
2020~\cite[Theorem~1.1]{GG20}. Hamaker, Pechenik, Speyer, and Weigandt
subsequently proved the determinant formula~\cite[Theorem~1.2]{HPSW20}.
Here the chains are obtained directly by
changing one coordinate of the forest code.

Our second main result describes when a lower interval is a Boolean lattice,
the inclusion order on subsets of a finite set. Tenner~\cite[Theorem~4.3]{Ten07}
proved that the principal ideal below a permutation $w$ in the strong Bruhat
order is Boolean if and only if $w$ avoids the patterns $321$ and $3412$.
Our classification gives a corresponding forbidden-configuration criterion
for forests.

\begin{thmB}[Boolean classification; Theorem~\ref{thm:boolean-classification}]
The lower interval $[\emptyforest,F]$ is a Boolean lattice if and only if $F$ is
an ``up-down forest": a disjoint union of trees in each of which the internal
vertices form two straight paths of left and right children descending from the
root. In that case
\[
  [\emptyforest,F]\cong B_{|\supp(F)|}.
\]
\end{thmB}

In an up-down forest, deleting support points determines the lower forest
uniquely, so the interval elements correspond to subsets of the support.
In the other direction, a change of direction along two successive internal
edges gives a \emph{bent corner}. Isolating this corner produces one of the
three-vertex non-Boolean intervals. The Boolean cube and the obstruction are
displayed in Figures~\ref{fig:interval-201} and~\ref{fig:interval-12}.

To study lower intervals, we use decorated insertion words and root lists.
Subword support containment makes it possible to delete selected letters and
restore the others through Monk covers. We also describe lower covers by
common bases and path-shaped spine quotients
(Theorem~\ref{thm:all-covers-zigzag}). The chain partition alone does not determine these lower covers.

The topological question is whether every principal lower interval is
shellable. We give an explicit recursive coatom ordering for the hook tops
$(n,1)$ (Theorem~\ref{thm:n1-shellable}). For the general question,
code-raising covers place every forest below a forest with no zero gaps in
its code. These are \emph{zigzag forests}, normalized to start at $1$.
Since shellability passes to intervals, they suffice for the conjecture.

\begin{conjC}[Shellability; Conjecture~\ref{conj:zigzag-shellable}]
Every zigzag lower interval $[\emptyforest,Z]$ is shellable; consequently every
lower interval $[\emptyforest,F]$ is shellable.
\end{conjC}

Section~\ref{sec:prelim} introduces the poset and its polynomial and code
models. Section~\ref{sec:sperner} proves the chain partition and strong
Sperner property, then counts the rank levels. Sections~\ref{sec:lower-covers}
and~\ref{sec:boolean} develop the insertion tools and Boolean
classification. Section~\ref{sec:shellability} treats hook shellability and
the general conjecture.

\section{Forest Polynomials and the Quasisymmetric Bruhat Order}
\label{sec:prelim}

We recall indexed forests and forest polynomials from Nadeau--Tewari
\cite[Sections~3.1--3.2]{NT23}, then use the forest Monk rule to define the order.

\subsection{Indexed Forests and Support}

\begin{definition}[Indexed forests]
An \emph{indexed tree} with support $[a,b]$ is a full binary tree whose leaves cover the numberline positions labeled $a, a+1, \dots, b + 1$. An \emph{indexed
forest} is a disjoint union of indexed trees with disjoint leaves.
\end{definition}

For a vertex $v$ in an indexed forest, let $\rho(v)$ be the leftmost number-line
label in the subtree rooted at $v$.

\begin{definition}[Support]
\label{def:support}
Draw an indexed forest above its number line. The \emph{support} $\supp(F)$ is
the set of number-line positions occupied by the component supports of $F$.
Thus a connected tree occupying the consecutive positions $a,\dots,b$ has
support $[a,b]$, and the support of a forest is the disjoint union of the
supports of its connected components.
\end{definition}

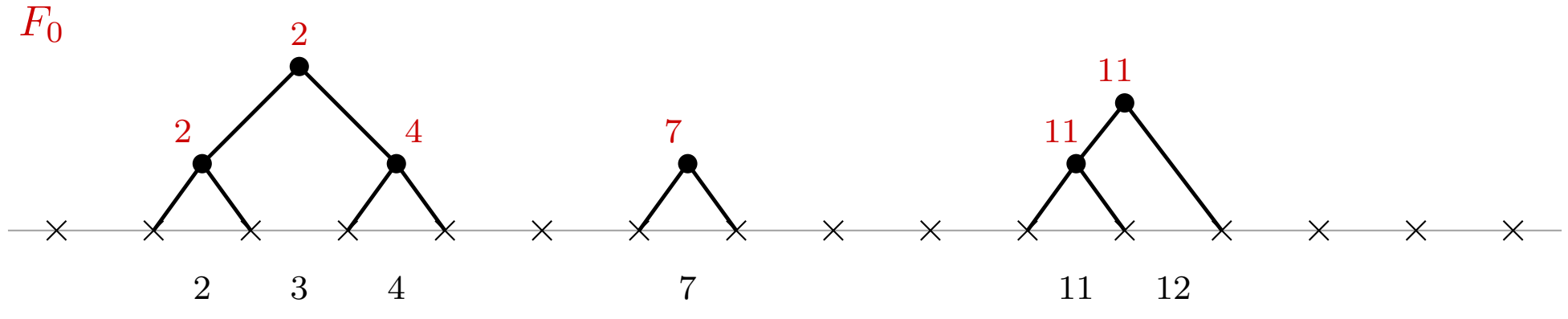
\begin{figure}[ht]
\centering
\begin{tikzpicture}[
  x=0.8cm,
  y=1.0cm,
  every node/.style={font=\small},
  forestedge/.style={line width=0.9pt},
  dot/.style={circle,fill=black,inner sep=1.6pt}
]
  \draw[gray!70] (0.5,0) -- (16.5,0);
  \foreach \i in {1,...,16} {
    \draw (\i-0.10,0.08) -- (\i+0.10,-0.08);
    \draw (\i-0.10,-0.08) -- (\i+0.10,0.08);
  }

  \node[draw=none] at (2.5,-0.47) {\scriptsize 2};
  \node[draw=none] at (3.5,-0.47) {\scriptsize 3};
  \node[draw=none] at (4.5,-0.47) {\scriptsize 4};

  \node[draw=none] at (7.5,-0.47) {\scriptsize 7};

  \node[draw=none] at (11.5,-0.47) {\scriptsize 11};
  \node[draw=none] at (12.5,-0.47) {\scriptsize 12};

  \coordinate (a) at (2.5,0.55);
  \coordinate (b) at (4.5,0.55);
  \coordinate (c) at (3.5,1.35);
  \draw[forestedge] (a) -- (2,0);
  \draw[forestedge] (a) -- (3,0);
  \draw[forestedge] (b) -- (4,0);
  \draw[forestedge] (b) -- (5,0);
  \draw[forestedge] (c) -- (a);
  \draw[forestedge] (c) -- (b);
  \node[dot] at (a) {};
  \node[dot] at (b) {};
  \node[dot] at (c) {};
  \node[draw=none,text=red!80!black] at (2.30,0.82) {\scriptsize 2};
  \node[draw=none,text=red!80!black] at (3.50,1.62) {\scriptsize 2};
  \node[draw=none,text=red!80!black] at (4.68,0.82) {\scriptsize 4};

  \coordinate (d) at (7.5,0.55);
  \draw[forestedge] (d) -- (7,0);
  \draw[forestedge] (d) -- (8,0);
  \node[dot] at (d) {};
  \node[draw=none,text=red!80!black] at (7.35,0.82) {\scriptsize 7};

  \coordinate (e) at (11.5,0.55);
  \coordinate (f) at (12.0,1.05);
  \draw[forestedge] (e) -- (11,0);
  \draw[forestedge] (e) -- (12,0);
  \draw[forestedge] (f) -- (e);
  \draw[forestedge] (f) -- (13,0);
  \node[dot] at (e) {};
  \node[dot] at (f) {};
  \node[draw=none,text=red!80!black] at (11.35,0.82) {\scriptsize 11};
  \node[draw=none,text=red!80!black] at (11.90,1.32) {\scriptsize 11};

  \node[draw=none,text=red!80!black] at (0.85,1.70) {$F_0$};
\end{tikzpicture}
\caption{An indexed forest $F_0$, adapted from~\cite[Figure~2]{NT23}.}
\label{fig:indexed-forest}
\end{figure}

In Figure~\ref{fig:indexed-forest}, the labels below the number line mark the
support positions, and the red labels record the leftmost descendant of each
internal vertex. For the left component, for example, the red labels are $2,2,4$.

\subsection{Forest Polynomials}

Forest polynomials are defined by summing over labelings of the internal
vertices.

\begin{definition}[Valid labelings and forest polynomials]
\label{def:forest-polynomial}
Let $F$ be an indexed forest, and let $\IN(F)$ denote its set of internal
vertices. A \emph{valid labeling} of $F$ is a map
\[
  \lambda:\IN(F)\to\ZZ_{>0}
\]
such that for every internal vertex $v$:
\begin{enumerate}[label=(\roman*)]
  \item $1\leq \lambda(v)\leq \rho(v)$;
  \item if $u$ is the left child of $v$, then $\lambda(v)\leq \lambda(u)$;
  \item if $u$ is the right child of $v$, then $\lambda(v)<\lambda(u)$.
\end{enumerate}
The \emph{forest polynomial} of $F$ is
\[
  P_F=\sum_{\lambda\text{ valid}}\prod_{v\in\IN(F)}x_{\lambda(v)}.
\]
\end{definition}

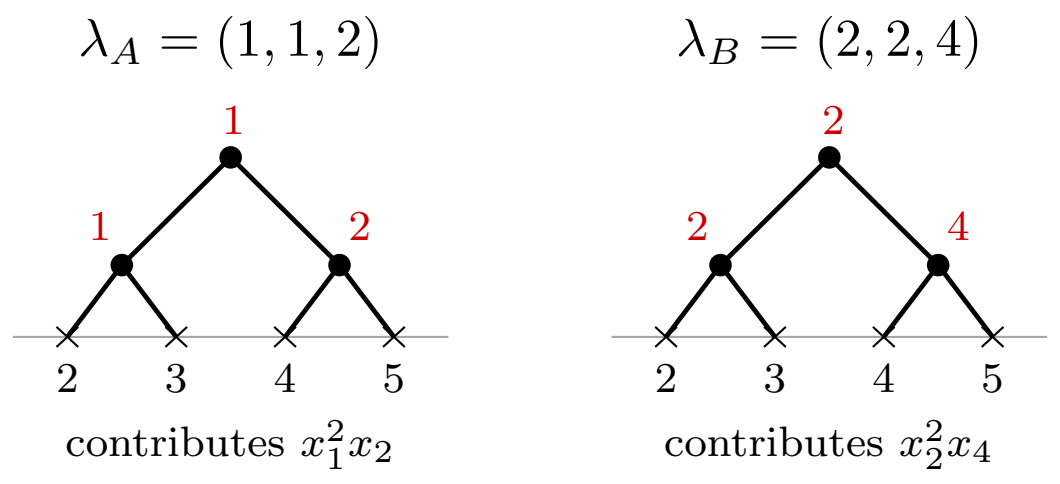
\begin{figure}[ht]
\centering
\begin{tikzpicture}[
  x=0.9cm,
  y=0.95cm,
  every node/.style={font=\small},
  treeedge/.style={line width=0.9pt},
  dot/.style={circle,fill=black,inner sep=1.5pt}
]
  \draw[gray!70] (0.0,-0.05) -- (3.2,-0.05);
  \foreach \x in {0.4,1.2,2.0,2.8} {
    \draw (\x-0.08,0.02) -- (\x+0.08,-0.12);
    \draw (\x-0.08,-0.12) -- (\x+0.08,0.02);
  }
  \node[draw=none] at (0.4,-0.34) {\scriptsize 2};
  \node[draw=none] at (1.2,-0.34) {\scriptsize 3};
  \node[draw=none] at (2.0,-0.34) {\scriptsize 4};
  \node[draw=none] at (2.8,-0.34) {\scriptsize 5};

  \draw[gray!70] (4.4,-0.05) -- (7.6,-0.05);
  \foreach \x in {4.8,5.6,6.4,7.2} {
    \draw (\x-0.08,0.02) -- (\x+0.08,-0.12);
    \draw (\x-0.08,-0.12) -- (\x+0.08,0.02);
  }
  \node[draw=none] at (4.8,-0.34) {\scriptsize 2};
  \node[draw=none] at (5.6,-0.34) {\scriptsize 3};
  \node[draw=none] at (6.4,-0.34) {\scriptsize 4};
  \node[draw=none] at (7.2,-0.34) {\scriptsize 5};

  \node at (1.6,2.0) {$\lambda_A=(1,1,2)$};
  \coordinate (A3) at (1.6,1.2);
  \coordinate (A2) at (0.8,0.45);
  \coordinate (A4) at (2.4,0.45);
  \draw[treeedge] (A3) -- (A2);
  \draw[treeedge] (A3) -- (A4);
  \draw[treeedge] (A2) -- (0.4,-0.05);
  \draw[treeedge] (A2) -- (1.2,-0.05);
  \draw[treeedge] (A4) -- (2.0,-0.05);
  \draw[treeedge] (A4) -- (2.8,-0.05);
  \node[dot] at (A2) {};
  \node[dot] at (A3) {};
  \node[dot] at (A4) {};
  \node[draw=none,text=red!80!black] at (0.64,0.72) {\scriptsize 1};
  \node[draw=none,text=red!80!black] at (1.62,1.46) {\scriptsize 1};
  \node[draw=none,text=red!80!black] at (2.55,0.72) {\scriptsize 2};
  \node at (1.6,-0.80) {\scriptsize contributes $x_1^2x_2$};

  \node at (6.0,2.0) {$\lambda_B=(2,2,4)$};
  \coordinate (B3) at (6.0,1.2);
  \coordinate (B2) at (5.2,0.45);
  \coordinate (B4) at (6.8,0.45);
  \draw[treeedge] (B3) -- (B2);
  \draw[treeedge] (B3) -- (B4);
  \draw[treeedge] (B2) -- (4.8,-0.05);
  \draw[treeedge] (B2) -- (5.6,-0.05);
  \draw[treeedge] (B4) -- (6.4,-0.05);
  \draw[treeedge] (B4) -- (7.2,-0.05);
  \node[dot] at (B2) {};
  \node[dot] at (B3) {};
  \node[dot] at (B4) {};
  \node[draw=none,text=red!80!black] at (5.03,0.72) {\scriptsize 2};
  \node[draw=none,text=red!80!black] at (6.03,1.46) {\scriptsize 2};
  \node[draw=none,text=red!80!black] at (6.95,0.72) {\scriptsize 4};
  \node at (6.0,-0.80) {\scriptsize contributes $x_2^2x_4$};
\end{tikzpicture}
\caption{Two valid labelings of $T_0$ and their monomials, adapted from~\cite[Figure~3, Example~3.5]{NT23}.}
\label{fig:two-labelings}
\end{figure}

Figure~\ref{fig:two-labelings} shows two valid labelings of the left component
$T_0$ of Figure~\ref{fig:indexed-forest}. In the displayed triples the entries
are read as (root, left child, right child).

The eight valid labelings of this three-vertex tree give
\begin{align*}
P_{T_0}
  ={}&x_1^2x_2+x_1^2x_3+x_1^2x_4+x_1x_2^2\\
    &+x_1x_2x_3+x_1x_2x_4+x_2^2x_3+x_2^2x_4.
\end{align*}
Thus the two labelings in Figure~\ref{fig:two-labelings} contribute the first
and last monomials in the expansion.

Disconnected components behave independently. If
$F=F_1\sqcup\cdots\sqcup F_r$ is the decomposition into connected components,
then
\[
  P_F=P_{F_1}\cdots P_{F_r}.
\]

\subsection{Codes}

\begin{definition}[Code]
\label{def:code}
The \emph{code} of a forest records the number of internal vertices with each
leftmost descendant:
\[
  c_i(F)=\#\{v\in\IN(F):\rho(v)=i\},
  \qquad
  \code(F)=(c_1,c_2,\dots).
\]
Trailing zeros are ignored. We write $F(\code)$ for the unique forest with
code $\code$~\cite[Proposition~3.3]{NT23}.
\end{definition}

For example, a tree with code $(2,0,1,1)$ has two vertices whose leftmost
descendant is $1$, no vertex whose leftmost descendant is $2$, and one vertex
over each of $3$ and $4$.

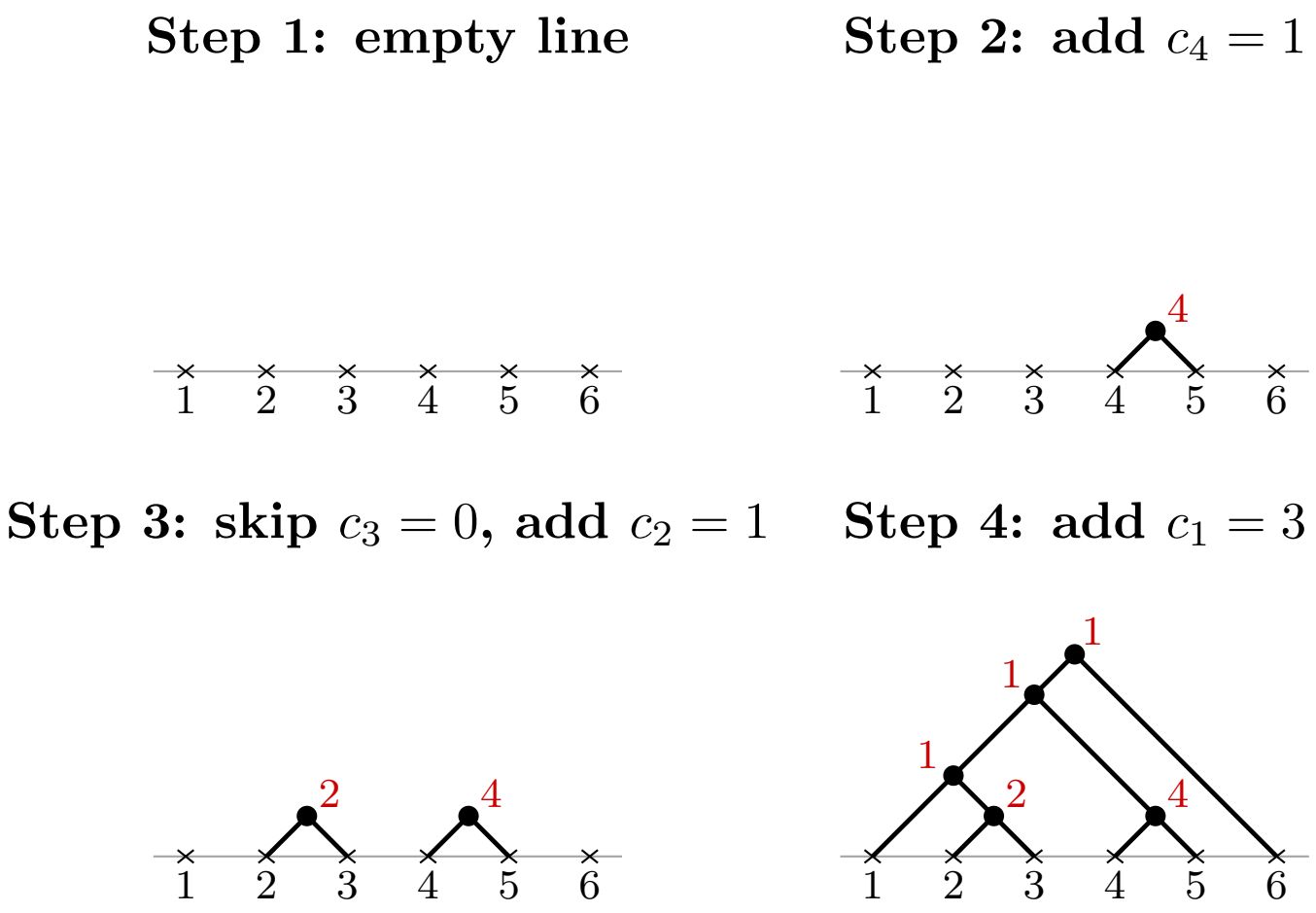
\begin{figure}[H]
\centering
\begin{tikzpicture}[
  x=0.55cm,
  y=0.55cm,
  every node/.style={font=\small},
  forestedge/.style={line width=0.9pt},
  dot/.style={circle,fill=black,inner sep=1.4pt}
]
  \begin{scope}[shift={(0,0)}]
    \node at (3.5,4.1) {\textbf{Step 1: empty line}};
    \draw[gray!70] (0.6,0) -- (6.4,0);
    \foreach \x/\lab in {1/1,2/2,3/3,4/4,5/5,6/6} {
      \draw (\x-0.10,0.08) -- (\x+0.10,-0.08);
      \draw (\x-0.10,-0.08) -- (\x+0.10,0.08);
      \node[draw=none] at (\x,-0.35) {\scriptsize \lab};
    }
  \end{scope}

  \begin{scope}[shift={(8.5,0)}]
    \node at (3.5,4.1) {\textbf{Step 2: add $c_4=1$}};
    \draw[gray!70] (0.6,0) -- (6.4,0);
    \foreach \x/\lab in {1/1,2/2,3/3,4/4,5/5,6/6} {
      \draw (\x-0.10,0.08) -- (\x+0.10,-0.08);
      \draw (\x-0.10,-0.08) -- (\x+0.10,0.08);
      \node[draw=none] at (\x,-0.35) {\scriptsize \lab};
    }
    \coordinate (a4) at (4.5,0.5);
    \draw[forestedge] (a4) -- (4,0);
    \draw[forestedge] (a4) -- (5,0);
    \node[dot] at (a4) {};
    \node[draw=none,text=red!80!black] at (4.78,0.78) {\scriptsize 4};
  \end{scope}

  \begin{scope}[shift={(0,-6.0)}]
    \node at (3.5,4.1) {\textbf{Step 3: skip $c_3=0$, add $c_2=1$}};
    \draw[gray!70] (0.6,0) -- (6.4,0);
    \foreach \x/\lab in {1/1,2/2,3/3,4/4,5/5,6/6} {
      \draw (\x-0.10,0.08) -- (\x+0.10,-0.08);
      \draw (\x-0.10,-0.08) -- (\x+0.10,0.08);
      \node[draw=none] at (\x,-0.35) {\scriptsize \lab};
    }
    \coordinate (b2) at (2.5,0.5);
    \coordinate (b4) at (4.5,0.5);
    \draw[forestedge] (b2) -- (2,0);
    \draw[forestedge] (b2) -- (3,0);
    \draw[forestedge] (b4) -- (4,0);
    \draw[forestedge] (b4) -- (5,0);
    \node[dot] at (b2) {};
    \node[dot] at (b4) {};
    \node[draw=none,text=red!80!black] at (2.78,0.78) {\scriptsize 2};
    \node[draw=none,text=red!80!black] at (4.78,0.78) {\scriptsize 4};
  \end{scope}

  \begin{scope}[shift={(8.5,-6.0)}]
    \node at (3.5,4.1) {\textbf{Step 4: add $c_1=3$}};
    \draw[gray!70] (0.6,0) -- (6.4,0);
    \foreach \x/\lab in {1/1,2/2,3/3,4/4,5/5,6/6} {
      \draw (\x-0.10,0.08) -- (\x+0.10,-0.08);
      \draw (\x-0.10,-0.08) -- (\x+0.10,0.08);
      \node[draw=none] at (\x,-0.35) {\scriptsize \lab};
    }
    \coordinate (c2) at (2.5,0.5);
    \coordinate (c4) at (4.5,0.5);
    \coordinate (c1) at (2.0,1.0);
    \coordinate (c3) at (3.0,2.0);
    \coordinate (c5) at (3.5,2.5);
    \draw[forestedge] (c5) -- (c3);
    \draw[forestedge] (c5) -- (6,0);
    \draw[forestedge] (c3) -- (c1);
    \draw[forestedge] (c3) -- (c4);
    \draw[forestedge] (c1) -- (1,0);
    \draw[forestedge] (c1) -- (c2);
    \draw[forestedge] (c2) -- (2,0);
    \draw[forestedge] (c2) -- (3,0);
    \draw[forestedge] (c4) -- (4,0);
    \draw[forestedge] (c4) -- (5,0);
    \node[dot] at (c1) {};
    \node[dot] at (c2) {};
    \node[dot] at (c3) {};
    \node[dot] at (c4) {};
    \node[dot] at (c5) {};
    \node[draw=none,text=red!80!black] at (1.68,1.26) {\scriptsize 1};
    \node[draw=none,text=red!80!black] at (2.78,0.78) {\scriptsize 2};
    \node[draw=none,text=red!80!black] at (2.72,2.26) {\scriptsize 1};
    \node[draw=none,text=red!80!black] at (4.78,0.78) {\scriptsize 4};
    \node[draw=none,text=red!80!black] at (3.72,2.78) {\scriptsize 1};
  \end{scope}
\end{tikzpicture}
\caption{Reconstructing a forest from the code $(3,1,0,1)$.}
\label{fig:code-to-tree}
\end{figure}

The factorization of Nadeau--Spink--Tewari~\cite[Proposition~4.2]{NST26}
reconstructs the forest from right to left. If $c_i=k$, the completed
forest has $k$ internal vertices with leftmost descendant $i$; zero entries are
passed over as gaps. Figure~\ref{fig:code-to-tree} illustrates the construction
for the code $(3,1,0,1)$.

The code is also visible inside the forest polynomial. Order monomials
lexicographically by their exponent vectors, read from left to right, with a
smaller exponent at the first differing coordinate considered smaller.

\begin{fact}[Lex-minimal monomial; {\cite[Proof of Proposition~3.6]{NT23}, \cite[Proposition~5.2]{NST26}}]
\label{fact:lex-min}
For every forest $F$, the lexicographically smallest monomial of $\PF$ is
$x^{\code(F)}$, and it appears with coefficient $1$. In particular, distinct
forest polynomials have distinct lex-smallest monomials.
\end{fact}

These distinct leading monomials immediately imply that the forest polynomials
$\{P_F\}$ form a basis of $\ZZ[x_1,x_2,\dots]$~\cite[Proposition~6.7]{NST26}.

\subsection{The Forest Monk Rule}

In the Schubert setting, Monk's formula describes multiplication by a divisor
class and identifies the resulting terms with Bruhat covers~\cite{Monk59}. The
forest analogue is the following theorem.

\begin{theorem}[Forest Monk rule; Nadeau--Spink--Tewari~{\cite[Theorem~7.4]{NST26}}]
\label{thm:forest-monk-rule}
For every indexed forest $F$ and every $i\geq 1$, the product
\[
  (x_1+\cdots+x_i)P_F
\]
expands as a linear combination of forest polynomials where all coefficients are $0$ or $1$.
\end{theorem}

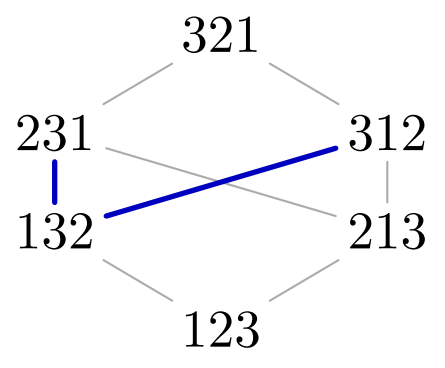
\begin{figure}[H]
\centering
\begin{tikzpicture}[x=1.1cm,y=.65cm,
  every node/.style={font=\small,inner sep=2pt},
  line cap=round]
  \node (e) at (0,0) {$123$};
  \node (a) at (-1,1) {$132$};
  \node (b) at (1,1) {$213$};
  \node (c) at (-1,2) {$231$};
  \node (d) at (1,2) {$312$};
  \node (t) at (0,3) {$321$};
  \draw[gray!65] (e)--(a) (e)--(b)
    (b)--(c) (b)--(d) (c)--(t) (d)--(t);
  \draw[blue!75!black,line width=1pt] (a)--(c) (a)--(d);
\end{tikzpicture}
\caption{The strong Bruhat order on $S_3$; $x_1\mathfrak S_{132}=\mathfrak S_{231}+\mathfrak S_{312}$ in blue.}
\label{fig:bruhat-monk-small}
\end{figure}

The classical model is the \emph{strong Bruhat order} on permutations. A cover
swaps two positions and increases the number of inversions by one, where
an inversion is a pair of entries in decreasing order~\cite{BB05}.
Writing $\mathfrak S_w$ for the Schubert polynomial of a permutation $w$,
the classical \emph{Monk rule} says that $(x_1+\cdots+x_i)\mathfrak S_w$
is the sum of the Schubert polynomials of its Bruhat covers that swap
positions $a\leq i<b$, each with coefficient one~\cite[Equation~(4.15$''$)]{Mac91}.
Figure~\ref{fig:bruhat-monk-small} shows a small example.

The forest Monk rule defines covers in the same way. The insertion model
in Section~\ref{sec:lower-covers} describes the individual covers.

\subsection{The Poset and Its Finite Windows}
\label{sec:first-structure}

\begin{definition}[Quasisymmetric Bruhat order]
\label{def:forest-monk-poset}
The \emph{quasisymmetric Bruhat order} is the poset whose elements are indexed forests
and whose cover relation is defined by
\[
  F\lessdot G
  \quad\Longleftrightarrow\quad
  \exists\,i\geq 1\ \text{such that}\
  P_G\text{ appears in }(x_1+\cdots+x_i)P_F.
\]
We write $E\leq F$ for the partial order generated by these covers.
\end{definition}

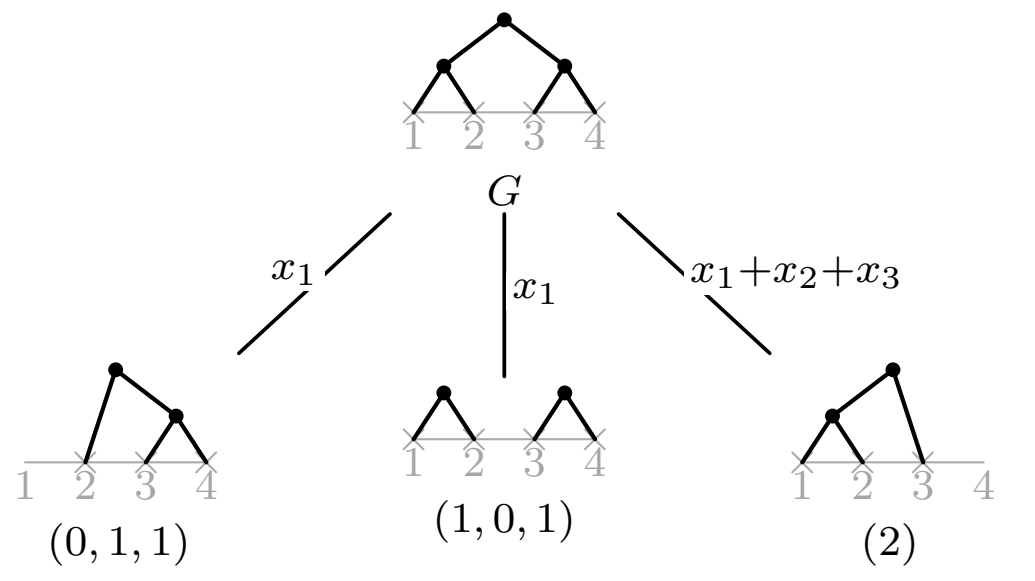
Figure~\ref{fig:first-covers} shows one forest $G$ and the forests it covers,
together with the Monk products that witness those covers.

\begin{figure}[H]
\centering
\begingroup
\newcommand{\covline}{%
  \draw[gray!65] (0,0)--(3,0);
  \foreach \x/\lab in {0/1,1/2,2/3,3/4}
    \node[below,gray!70] at (\x,0) {\scriptsize \lab};}
\newcommand{\covleaf}[1]{\node[gray!70] at (#1,0) {\scriptsize $\times$};}
\newcommand{\fCoverTop}{%
\begin{tikzpicture}[x=.40cm,y=.36cm,baseline=-0.6ex,line cap=round,line join=round]
  \covline
  \covleaf{0}\covleaf{1}\covleaf{2}\covleaf{3}
  \draw[thick] (1.5,1.7)--(.5,.85) (1.5,1.7)--(2.5,.85)
               (.5,.85)--(0,0) (.5,.85)--(1,0)
               (2.5,.85)--(2,0) (2.5,.85)--(3,0);
  \fill (1.5,1.7) circle (1.4pt);
  \fill (.5,.85) circle (1.4pt);
  \fill (2.5,.85) circle (1.4pt);
\end{tikzpicture}}
\newcommand{\fCoverLeft}{%
\begin{tikzpicture}[x=.40cm,y=.36cm,baseline=-0.6ex,line cap=round,line join=round]
  \covline
  \covleaf{1}\covleaf{2}\covleaf{3}
  \draw[thick] (1.5,1.7)--(1,0) (1.5,1.7)--(2.5,.85)
               (2.5,.85)--(2,0) (2.5,.85)--(3,0);
  \fill (1.5,1.7) circle (1.4pt);
  \fill (2.5,.85) circle (1.4pt);
\end{tikzpicture}}
\newcommand{\fCoverMid}{%
\begin{tikzpicture}[x=.40cm,y=.36cm,baseline=-0.6ex,line cap=round,line join=round]
  \covline
  \covleaf{0}\covleaf{1}\covleaf{2}\covleaf{3}
  \draw[thick] (.5,.85)--(0,0) (.5,.85)--(1,0)
               (2.5,.85)--(2,0) (2.5,.85)--(3,0);
  \fill (.5,.85) circle (1.4pt);
  \fill (2.5,.85) circle (1.4pt);
\end{tikzpicture}}
\newcommand{\fCoverRight}{%
\begin{tikzpicture}[x=.40cm,y=.36cm,baseline=-0.6ex,line cap=round,line join=round]
  \covline
  \covleaf{0}\covleaf{1}\covleaf{2}
  \draw[thick] (1.5,1.7)--(.5,.85) (1.5,1.7)--(2,0)
               (.5,.85)--(0,0) (.5,.85)--(1,0);
  \fill (1.5,1.7) circle (1.4pt);
  \fill (.5,.85) circle (1.4pt);
\end{tikzpicture}}
\begin{tikzpicture}[
  x=1cm, y=1cm,
  line cap=round, line join=round,
  forestnode/.style={align=center,inner sep=1.5pt},
  cover/.style={line width=.7pt},
  edgelabel/.style={fill=white,inner sep=1.2pt,font=\scriptsize}]
  \node[forestnode] (top) at (0,2.35) {\fCoverTop\\[1pt]{\scriptsize $G$}};
  \node[forestnode] (L) at (-2.55,0) {\fCoverLeft\\[1pt]{\scriptsize $(0,1,1)$}};
  \node[forestnode] (M) at (0,0) {\fCoverMid\\[1pt]{\scriptsize $(1,0,1)$}};
  \node[forestnode] (R) at (2.55,0) {\fCoverRight\\[1pt]{\scriptsize $(2)$}};
  \draw[cover] (L)--node[pos=.58,edgelabel,left] {$x_1$} (top);
  \draw[cover] (M)--node[pos=.52,edgelabel,right] {$x_1$} (top);
  \draw[cover] (R)--node[pos=.58,edgelabel,right] {$x_1{+}x_2{+}x_3$} (top);
  \node[anchor=west,align=left,font=\small] at (4.15,1.15)
    {$P_G=x_1^2x_2+x_1^2x_3$\\[3pt]
     $P_{(0,1,1)}=x_1x_2+x_1x_3+x_2x_3$\\[3pt]
     $P_{(1,0,1)}=x_1^2+x_1x_2+x_1x_3$\\[3pt]
     $P_{(2)}=x_1^2$\\[7pt]
     $x_1\,P_{(0,1,1)}=P_G+P_{(1,1,1)}$\\[3pt]
     $x_1\,P_{(1,0,1)}=P_G+P_{(3)}$\\[3pt]
     $(x_1{+}x_2{+}x_3)\,P_{(2)}=P_G+P_{(3)}$};
\end{tikzpicture}
\endgroup
\caption{The lower covers of $(2,0,1)$.}
\label{fig:first-covers}
\end{figure}

For a forest $F$, its \emph{principal lower interval} is
$[\emptyforest,F]=\{E:E\leq F\}$. For a finite set $S$ of positive integers,
let $\mathcal P_S$ be the induced subposet of forests supported in $S$.
Section~\ref{subsec:component-covers} gives the product decomposition for a
general finite support set.

Every forest polynomial is homogeneous of degree $|\IN(F)|$, so a Monk
cover increases this degree by one. The rank is therefore
\[
  \rk(F)=|\IN(F)|=|\supp(F)|=\deg P_F.
\]

\begin{definition}[Standard window]
\label{def:standard-window}
For the rank and Sperner results we work in the finite window
\[
  \mathcal{P}_n=\mathcal{P}_{[1,n]},
\]
the induced subposet of forests whose support is contained in $[1,n]$, graded
by rank $0,1,\dots,n$. Equivalently, $\mathcal{P}_n$ is the window of forests
whose leaves lie in $[1,n+1]$: a maximal element has full support $[1,n]$, hence
$n$ internal vertices and $n+1$ leaves filling the positions $1,\dots,n+1$.
\end{definition}

We write $\mathcal{P}_n^{(k)}$ for the rank-$k$ level of $\mathcal{P}_n$, the set
of its elements of rank $k$.

\section{Chains and the Strong Sperner Property}
\label{sec:sperner}

Raising the first code entry partitions $\mathcal P_n$ into saturated chains
ending in top rank. An antichain can meet each chain only once; more generally,
the partition gives the sharp bound for a union of antichains. We first
construct the chains, then evaluate the rank numbers.

\subsection{Code-Raising Covers}

Raising any single code entry produces a cover. Here $e_i$ denotes the $i$th standard unit vector, so $\code+e_i$
increments the $i$th code entry by one.

\begin{lemma}[Code-raising covers]
\label{lem:code-raising}
For every forest $F$ with $\code(F)=c$ and every $i\geq 1$,
\[
  F(c)\;\lessdot\;F(c+e_i).
\]
\end{lemma}

\begin{proof}
The least monomial identifies the required summand without expanding the
whole product. Since lexicographic order is preserved by multiplication,
Fact~\ref{fact:lex-min} gives
\[
  \min_{\mathrm{lex}}\bigl((x_1+\cdots+x_i)\PF\bigr)=x^{c+e_i}.
\]
The positive Monk expansion contains a summand $P_G$ with this monomial.
Its least monomial cannot be smaller, since it would then occur in the
product as well. Thus the least monomial of $P_G$ is $x^{c+e_i}$, and
Fact~\ref{fact:lex-min} identifies $\code(G)=c+e_i$. This gives the cover.
\end{proof}

\subsection{The Chain Partition}

Membership in $\mathcal{P}_n$ is characterized by bounds on the code's
suffix sums.

\begin{lemma}[Window code condition]
\label{lem:window-codes}
A code $\code=(c_1,c_2,\dots)$ has $F(\code)\in\mathcal{P}_n$ if and only if its
suffix sums satisfy
\[
  \sum_{j\geq i}c_j\leq n+1-i
  \qquad\text{for every }1\leq i\leq n+1.
\]
\end{lemma}

\begin{proof}
The suffix sum $g(i)=\sum_{j\geq i}c_j$ counts vertices whose leftmost
descendant leaf has label at least $i$. Each occupies a distinct inorder
support position at least $i$. If $\supp(F)\subseteq[1,n]$, there are only
$n+1-i$ such positions, giving the required bound.

Conversely, suppose $M=\max\supp(F)>n$. The rightmost connected component
has support $[a,M]$, so its $M-a+1$ vertices all have $\rho\geq a$.
If $a\leq n+1$, then
\[
  g(a)\geq M-a+1>n+1-a,
\]
contrary to the bound at $i=a$. If $a>n+1$, these vertices instead give
$g(n+1)>0$, contrary to the bound at $i=n+1$. Hence the suffix bounds force
$\supp(F)\subseteq[1,n]$.
\end{proof}

The \emph{first-entry map} is
\label{def:pumping-map}
\[
  \Phi(F(\code))=F(\code+e_1).
\]
By Lemma~\ref{lem:code-raising}, $F(\code)\lessdot\Phi(F(\code))$.
Its support behavior keeps these covers inside the finite window.

Write $\mex(S)$ for the least positive integer not in a set $S$. The support
point added by $\Phi$ is always this leftmost vacancy.

\begin{lemma}[Mex law]
\label{lem:mex}
For every code $c$,
\[
  \supp\bigl(F(c+e_1)\bigr)=\supp\bigl(F(c)\bigr)\cup\{\mex(\supp(F(c)))\}.
\]
Consequently, if $F(c)\in\mathcal{P}_n$ has rank below $n$, then
$F(c+e_1)\in\mathcal{P}_n$; below the top rank, $\Phi$ maps $\mathcal{P}_n$ into
itself.
\end{lemma}

\begin{proof}
The codes $c$ and $c+e_1$ agree for every $i\geq 2$, so their right-to-left
reconstructions~\cite[Proposition~4.2]{NST26} coincide until the extra
installation at position $1$. That step creates exactly the leftmost vacant support point of $F(c)$. If
$F(c)$ has rank below $n$, some point of $[1,n]$ is missing, so the new point
lies in the window and $F(c+e_1)$ stays in $\mathcal{P}_n$.
\end{proof}
\begin{observation}[Multiple maximal elements]
\label{obs:multiple-maxima}
The maximal elements of $\mathcal P_n$ are exactly the forests with full
support $[1,n]$. There are $\cat_n$ such maximal elements.
\end{observation}

\begin{proof}
A full-support forest has rank $n$, so no forest in $\mathcal P_n$ can cover
it. Conversely, if a forest in $\mathcal P_n$ has a missing
support point in $[1,n]$, Lemma~\ref{lem:mex} gives an upper cover adding
one missing point, so the forest is not maximal. The full-support forests are
the full binary tree shapes on $n$ internal vertices, counted by the Catalan
number $\cat_n$.
\end{proof}

\begin{figure}[H]
\centering
\begin{tikzpicture}[
    x=.95cm,y=1.16cm,
    cover/.style={gray!35,line width=.45pt},
    highway/.style={blue!70!black,line width=1.45pt},
    entry/.style={fill=white,inner sep=1.5pt,font=\small},
    ranklabel/.style={anchor=east,font=\small}]
  \coordinate (z) at (0,0);

  \coordinate (u001) at (-4,1);
  \coordinate (u01)  at (0,1);
  \coordinate (u1)   at (4,1);

  \coordinate (v011) at (-4,2);
  \coordinate (v02)  at (-2,2);
  \coordinate (v101) at (0,2);
  \coordinate (v11)  at (2,2);
  \coordinate (v2)   at (4,2);

  \coordinate (w111) at (-4,3);
  \coordinate (w12)  at (-2,3);
  \coordinate (w201) at (0,3);
  \coordinate (w21)  at (2,3);
  \coordinate (w3)   at (4,3);

  \draw[cover]
    (z)--(u001) (z)--(u1) (z)--(u01)
    (u001)--(v011) (u001)--(v02) (u001)--(v101)
    (u01)--(v011) (u01)--(v02) (u01)--(v11) (u01)--(v2)
    (u1)--(v101) (u1)--(v11) (u1)--(v2)
    (v011)--(w111) (v011)--(w12) (v011)--(w201)
    (v02)--(w12) (v02)--(w21) (v02)--(w3)
    (v101)--(w111) (v101)--(w12) (v101)--(w201)
    (v101)--(w21) (v101)--(w3)
    (v11)--(w111) (v11)--(w12) (v11)--(w21)
    (v2)--(w201) (v2)--(w21) (v2)--(w3);

  \draw[highway]
    (z)--(u1)--(v2)--(w3)
    (u001)--(v101)--(w201)
    (u01)--(v11)--(w21)
    (v011)--(w111)
    (v02)--(w12);

  \node[ranklabel] at (-5.05,0) {rank $0$};
  \node[ranklabel] at (-5.05,1) {rank $1$};
  \node[ranklabel] at (-5.05,2) {rank $2$};
  \node[ranklabel] at (-5.05,3) {rank $3$};

  \node[entry] at (z) {$\emptyforest$};
  \node[entry] at (u001) {$(0,0,1)$};
  \node[entry] at (u1) {$(1)$};
  \node[entry] at (u01) {$(0,1)$};
  \node[entry] at (v101) {$(1,0,1)$};
  \node[entry] at (v011) {$(0,1,1)$};
  \node[entry] at (v2) {$(2)$};
  \node[entry] at (v02) {$(0,2)$};
  \node[entry] at (v11) {$(1,1)$};
  \node[entry] at (w201) {$(2,0,1)$};
  \node[entry] at (w111) {$(1,1,1)$};
  \node[entry] at (w3) {$(3)$};
  \node[entry] at (w12) {$(1,2)$};
  \node[entry] at (w21) {$(2,1)$};
\end{tikzpicture}
\caption{The pumping-chain partition of $\mathcal{P}_3$ in blue.}
\label{fig:pumping-highways}
\end{figure}
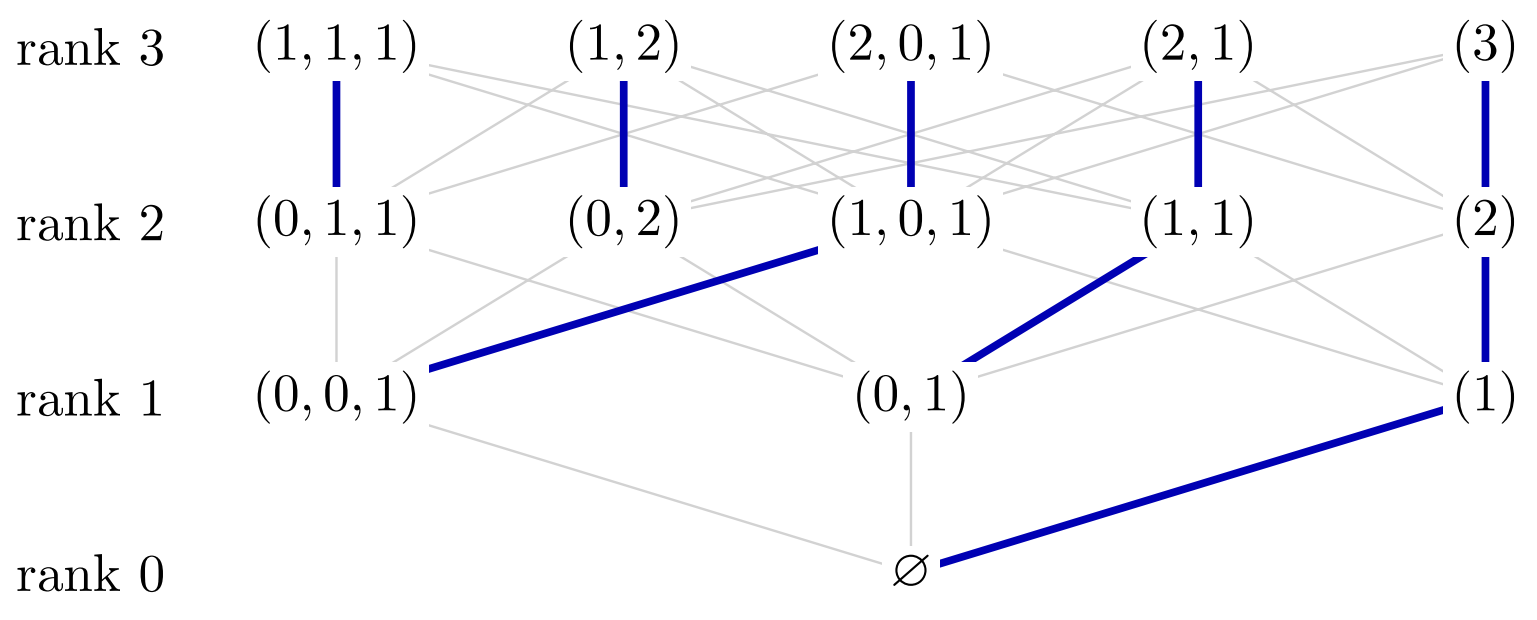

Two codes are \emph{tail equivalent} when they agree except possibly in their
first entry:
\label{def:tail-equivalence}
\[
  (c_1,c_2,c_3,\dots)\sim(d_1,c_2,c_3,\dots).
\]
Each tail-equivalence class is one pumping chain. In $\mathcal{P}_3$ the five
classes are the blue chains of Figure~\ref{fig:pumping-highways}: every forest
lies on exactly one, and each chain is saturated and ends at rank $3$. For
instance,
\[
  F(0,0,1)\;\lessdot\;F(1,0,1)\;\lessdot\;F(2,0,1).
\]

\begin{samepage}
\begin{theorem}[Pumping chain partition]
\label{thm:chain-partition}
The first-entry map $\Phi$ has the following properties.
\begin{enumerate}[label=(\roman*)]
  \item $\Phi$ is defined on every element of $\mathcal{P}_n$ of rank below
    $n$, maps into $\mathcal{P}_n$, and $F\lessdot\Phi(F)$.
  \item $\Phi$ is injective.
  \item Iterating $\Phi$ partitions $\mathcal{P}_n$ into saturated chains
  \[
    F(0,c_2,c_3,\dots)\;\lessdot\;F(1,c_2,c_3,\dots)\;\lessdot\;
    F(2,c_2,c_3,\dots)\;\lessdot\;\cdots,
  \]
  indexed by the codes with $c_1=0$.
  \item Every chain terminates in the top rank level $\mathcal{P}_n^{(n)}$.
\end{enumerate}
\end{theorem}
\end{samepage}
\begin{proof}
Fix a tail $(c_2,c_3,\dots)$ occurring in the window, and set
$s=\sum_{j\geq2}c_j$. In Lemma~\ref{lem:window-codes}, only the first
suffix bound involves $c_1$, so its allowable values are exactly
\[
  c_1=0,1,\ldots,n-s.
\]
Raising $c_1$ walks through these values by covers, by
Lemma~\ref{lem:code-raising}, and the last forest has rank $(n-s)+s=n$.
Different tails give disjoint chains, and every code belongs to one of
them. Finally, $c+e_1$ determines $c$, so $\Phi$ is injective. These
observations give all four assertions.
\end{proof}

\subsection{Sperner Consequences}

\begin{definition}[Strong Sperner property]
\label{def:strong-sperner}
A ranked finite poset is \emph{strongly Sperner} if, for every $r\geq 1$, the
largest union of $r$ antichains has size equal to the sum of the $r$ largest rank
levels; see Engel~\cite{Eng97} for the general theory.
If $r$ exceeds the number of rank levels, the missing levels are taken to be empty.
\end{definition}

The chain partition bounds the size of each union of antichains. In $\mathcal{P}_3$, the five
blue chains of Figure~\ref{fig:pumping-highways} already force every antichain
to have size at most five, and every union of two antichains to have size at
most ten.

\begin{corollary}[Strong Sperner property]
\label{cor:strong-sperner}
The standard window $\mathcal{P}_n$ is strongly Sperner, and its largest antichains
have size $\cat_{n}$.
\end{corollary}

\begin{proof}
The rank levels weakly increase because $\Phi$ injects each level into
the next. Let $A$ be a union of $r$ antichains. In a pumping chain $C$,
it can occupy at most $\min\{r,|C|\}$ vertices. Since every chain ends
at rank $n$, this is exactly the number of its vertices in the top $r$
levels. Summing over the chain partition gives
\[
  |A|\leq\sum_C\min\{r,|C|\}
  =\sum_{j=\max\{0,n-r+1\}}^n |\mathcal P_n^{(j)}|.
\]
The top $r$ levels attain this bound and are the $r$ largest levels.
Taking $r=1$ gives the largest-antichain size $\cat_n$.
\end{proof}

\begin{remark}
The chain through any forest is found by reducing $c_1$ to $0$, then raising it
until the support fills, as in the blue edges of
Figure~\ref{fig:pumping-highways}. 
\end{remark}

\subsection{Ballot Rank Numbers}
\label{sec:rank-generating}

We recover the known enumeration of indexed forests by size~\cite{NT23},
recorded explicitly in~\cite[Remark~9.4]{BGNST26}. These ballot numbers
also give the graded dimensions of the quasisymmetric coinvariant ring
in $n+1$ variables~\cite[Theorem~5.1]{ABB04}. We use the cycle lemma of
Dvoretzky and Motzkin~\cite{DM47}.

\begin{lemma}[Cycle lemma]
\label{lem:cycle}
Let $b=(b_1,\dots,b_N)$ be a sequence of integers, each at most $1$, with sum
$s\geq1$. Then exactly $s$ of its $N$ indexed cyclic rotations have all
partial sums positive.
\end{lemma}

This form follows from the usual $+1,-1$ cycle lemma~\cite[Section~1]{DZ90}
by replacing each $b_i$ with one $+1$ followed by $1-b_i$ copies of $-1$.

\begin{theorem}[Rank numbers]
\label{thm:rank-numbers}
For $0\leq k\leq n$,
\[
  |\mathcal{P}_n^{(k)}|=\frac{n+1-k}{n+1}\binom{n+k}{k}.
\]
Consequently the rank generating function is
\[
  R_n(q)=\sum_{k=0}^{n}\frac{n+1-k}{n+1}\binom{n+k}{k}\,q^k.
\]
\end{theorem}

\begin{proof}
By Lemma~\ref{lem:window-codes}, a rank-$k$ code has $c_j=0$ for $j\geq n+1$.
Reversing its remaining entries, $d_p=c_{n+1-p}$, identifies
$\mathcal P_n^{(k)}$ with the nonnegative sequences $(d_1,\dots,d_n)$ of
sum $k$ satisfying
\[
  d_1+\cdots+d_p\leq p\qquad(1\leq p\leq n).
\]
Set $N=n+1$ and $s=N-k$. The encoding
\[
  b=(1,\;1-d_1,\;\dots,\;1-d_n)
\]
has entries at most $1$, sum $s$, and positive partial sums exactly when
these prefix bounds hold. Conversely, any integer sequence with entries
at most $1$ and positive partial sums must start with $1$. Thus this is a
bijection with the length-$N$ sequences of sum $s$ having those properties.

Without the positivity condition, writing $b_i=1-e_i$ identifies these
sequences with weak compositions of $k$ into $N$ parts, of which there are
$\binom{n+k}{k}$. This set is closed under rotation. Count pairs $(b,j)$,
where $1\leq j\leq N$ and the rotation starting at $j$ has positive partial
sums. Lemma~\ref{lem:cycle} gives $s$ such indexed rotations for each $b$.
Each positive sequence is obtained from exactly one pair for each fixed
$j$, even when some rotations coincide. Therefore
\[
  (n+1)|\mathcal P_n^{(k)}|=(n+1-k)\binom{n+k}{k},
\]
which gives the rank formula.
\end{proof}

Table~\ref{tab:rank-levels} lists several rows of the Catalan ballot triangle.

\begin{table}[ht]
\centering
\begin{tabular}{cl}
\toprule
$n$ & rank level sizes $|\mathcal{P}_n^{(0)}|,\dots,|\mathcal{P}_n^{(n)}|$ \\
\midrule
$4$ & $1,4,9,14,14$ \\
$5$ & $1,5,14,28,42,42$ \\
$6$ & $1,6,20,48,90,132,132$ \\
$7$ & $1,7,27,75,165,297,429,429$ \\
$8$ & $1,8,35,110,275,572,1001,1430,1430$ \\
\bottomrule
\end{tabular}
\caption{The rank sizes of $\mathcal{P}_n$: the Catalan ballot triangle.}
\label{tab:rank-levels}
\end{table}

The total count below was also obtained by the cycle lemma
in~\cite[Section~3.1]{NST26}, where $\mathcal P_n$ is denoted
$\mathsf{For}_{n+1}$.

\begin{corollary}[Catalan sizes]
\label{cor:catalan-size}
The standard window has cardinalities
\[
  |\mathcal{P}_n|=\cat_{n+1},
  \qquad
  |\mathcal{P}_n^{(n)}|=\cat_n.
\]
\end{corollary}

\begin{proof}
Set $k=n$ in Theorem~\ref{thm:rank-numbers}, and sum its rank formula using
the classical ballot row-sum identity.
\end{proof}

\begin{corollary}[Unimodality]
\label{cor:unimodality}
For $n\ge1$, the rank numbers of $\mathcal{P}_n$ strictly increase up to
the last two levels, which are equal:
\[
  |\mathcal{P}_n^{(0)}|<|\mathcal{P}_n^{(1)}|<\cdots
  <|\mathcal{P}_n^{(n-1)}|=|\mathcal{P}_n^{(n)}|.
\]
\end{corollary}

\begin{remark}[Number of chains]
For $n\ge1$, the pumping chains are indexed by forests with $c_1=0$,
equivalently with $1\notin\supp(F)$. Shifting these forests down by one
identifies them with $\mathcal P_{n-1}$, so the number of chains is
\[
  |\mathcal P_{n-1}|=\cat_n=|\mathcal P_n^{(n)}|.
\]
For $n=0$, the empty forest forms the unique chain.
\end{remark}

\section{Lower Intervals through Insertion}
\label{sec:lower-covers}

Our next goal is the Boolean classification in Theorem~B. For an up-down
forest $F$, we will show that every subset of $\supp(F)$ determines
exactly one forest below $F$, with covers corresponding to deleting a
single support point. Conversely, for a forest that is not up-down, we
will isolate a bent corner as a lower element whose own lower interval
is not Boolean. Such an interval cannot occur inside a Boolean lattice.

Both arguments require a concrete way to move downward in the order.
In this section we represent forests by insertion words and study what happens when
letters are deleted. Component decomposition will reduce the problem
to one tree at a time; support-removing covers will produce the subsets
needed for the Boolean model; and subword support containment will
allow us to isolate the forbidden configuration. These tools are then used in Section~\ref{sec:boolean}.

\subsection{Decorated Insertion Words}

The decorated insertion of Nadeau and Tewari~\cite[Section~5.1 and Remark~5.2]{NT23}
generalizes the Sylvester correspondence between words and binary trees.
It reads a word one letter at a time and builds a forest.

A \emph{decorated letter} has the form $i^{[e]}$, where $i$ is its value and
$e$ breaks ties among letters of the same value. The letters are totally
ordered by value and then by decoration. Given a word $w$ in distinct decorated
letters, write $\OmegaF(w)$ for its insertion forest, denoted $F(w)$
in~\cite{NT23}; their $\Omega(w)$ denotes the parked set. Each new letter
becomes the root of the component it creates or enlarges.

\begin{figure}[H]
\centering
\begingroup
\newcommand{\leafline}[2]{%
  \draw[gray!65] (#1,0)--(#2,0);
  \foreach \x in {#1,...,#2} {
    \node[gray!70] at (\x,0) {\scriptsize $\times$};
    \node[below,gray!70] at (\x,0) {\scriptsize \x};
  }
}
\begin{tikzpicture}[
  x=.72cm,
  y=.72cm,
  every node/.style={font=\scriptsize},
  line cap=round,
  line join=round,
  forestedge/.style={line width=.85pt},
  newedge/.style={line width=.95pt,blue!75!black},
  dot/.style={circle,fill=black,inner sep=1.4pt},
  newdot/.style={circle,fill=blue!75!black,inner sep=1.6pt},
  title/.style={font=\scriptsize,align=center}
]
  \begin{scope}[shift={(0,0)}]
    \node[title] at (2.5,2.55) {empty spot};
    \leafline{1}{4}
    \coordinate (u) at (2.5,0.5);
    \draw[newedge] (u)--(2,0);
    \draw[newedge] (u)--(3,0);
    \node[newdot] at (u) {};
    \node[blue!75!black] at (2.5,0.95) {$2^{[1]}$};
  \end{scope}

  \begin{scope}[shift={(7.4,0)}]
    \node[title] at (2.5,2.55) {bridge a gap};
    \leafline{1}{4}
    \coordinate (lft) at (1.5,0.5);
    \coordinate (rgt) at (3.5,0.5);
    \coordinate (u)   at (2.5,1.5);
    \draw[forestedge] (lft)--(1,0);
    \draw[forestedge] (lft)--(2,0);
    \draw[forestedge] (rgt)--(3,0);
    \draw[forestedge] (rgt)--(4,0);
    \draw[newedge] (u)--(lft);
    \draw[newedge] (u)--(rgt);
    \node[dot] at (lft) {};
    \node[dot] at (rgt) {};
    \node[newdot] at (u) {};
    \node at (1.06,0.62) {$1^{[1]}$};
    \node at (3.94,0.62) {$3^{[1]}$};
    \node[blue!75!black] at (2.5,1.95) {$2^{[1]}$};
  \end{scope}

  \begin{scope}[shift={(0,-5.0)}]
    \node[title] at (2.5,2.95) {grow right};
    \leafline{1}{4}
    \coordinate (subn)    at (2.5,0.5);
    \coordinate (oldroot) at (2.0,1.0);
    \coordinate (u)       at (2.5,1.5);
    \draw[forestedge] (oldroot)--(1,0);
    \draw[forestedge] (oldroot)--(subn);
    \draw[forestedge] (subn)--(2,0);
    \draw[forestedge] (subn)--(3,0);
    \draw[newedge] (u)--(oldroot);
    \draw[newedge] (u)--(4,0);
    \node[dot] at (oldroot) {};
    \node[dot] at (subn) {};
    \node[newdot] at (u) {};
    \node at (1.46,1.16) {$1^{[1]}$};
    \node at (2.96,0.52) {$2^{[1]}$};
    \node[blue!75!black] at (2.5,1.95) {$2^{[2]}$};
  \end{scope}

  \begin{scope}[shift={(7.4,-5.0)}]
    \node[title] at (2.5,2.95) {grow left};
    \leafline{1}{4}
    \coordinate (subn)    at (2.5,0.5);
    \coordinate (oldroot) at (3.0,1.0);
    \coordinate (u)       at (2.5,1.5);
    \draw[forestedge] (oldroot)--(subn);
    \draw[forestedge] (oldroot)--(4,0);
    \draw[forestedge] (subn)--(2,0);
    \draw[forestedge] (subn)--(3,0);
    \draw[newedge] (u)--(1,0);
    \draw[newedge] (u)--(oldroot);
    \node[dot] at (oldroot) {};
    \node[dot] at (subn) {};
    \node[newdot] at (u) {};
    \node at (3.54,1.16) {$2^{[2]}$};
    \node at (2.04,0.52) {$2^{[1]}$};
    \node[blue!75!black] at (2.5,1.95) {$1^{[1]}$};
  \end{scope}
\end{tikzpicture}
\caption{The four local moves in decorated insertion.}
\label{fig:omega-local-moves}
\endgroup
\end{figure}
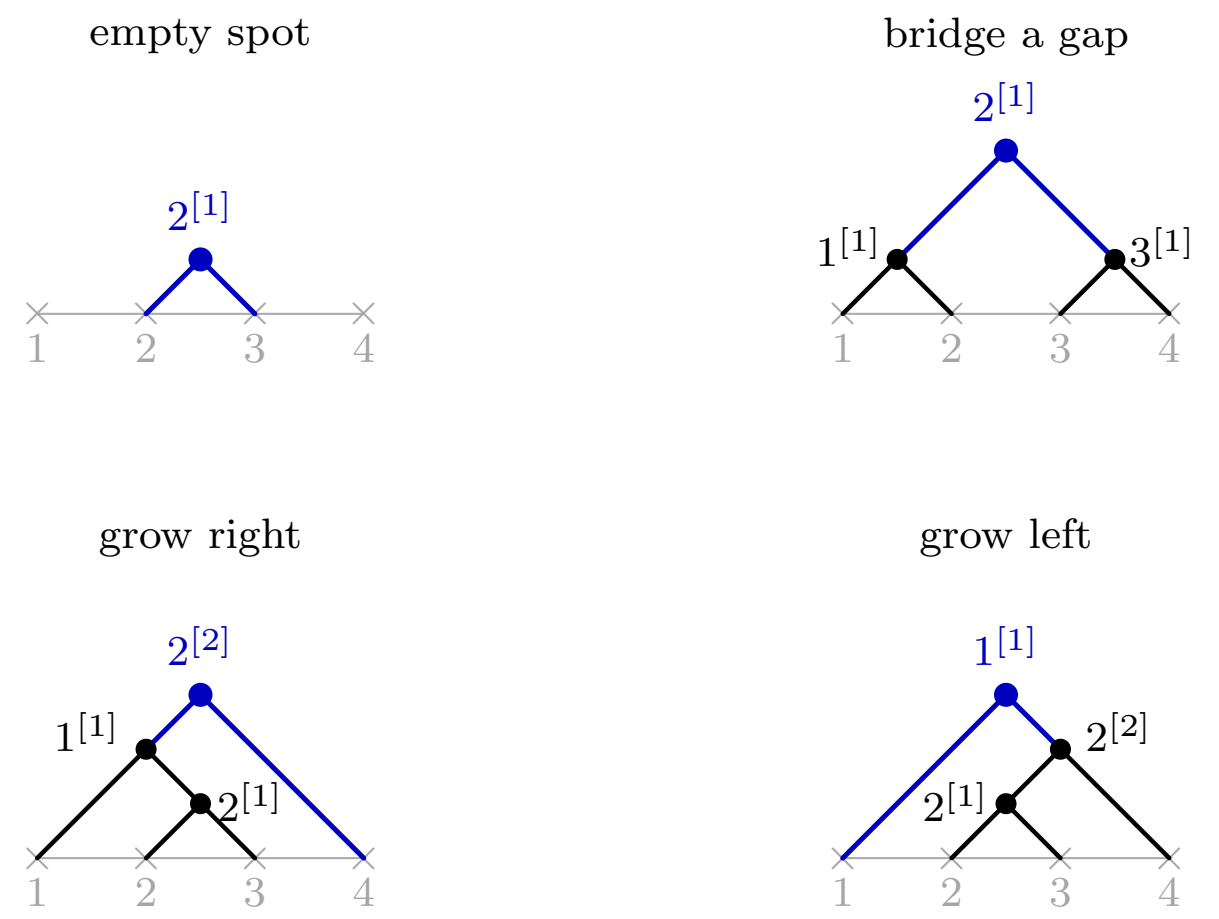

The four local moves are shown in Figure~\ref{fig:omega-local-moves}. Black
nodes belong to the existing forest, while the blue node and edges are created
by the newly inserted decorated letter.

Each insertion adds one support point, as follows.

\begin{fact}[Support growth under one insertion]
\label{fact:support-growth}
If $G$ is obtained from $F$ by inserting one decorated letter, then $G$ has
exactly one more internal vertex than $F$, and there is a unique support point
$x$ such that
\[
  \supp(G)=\supp(F)\sqcup\{x\}.
\]
If the inserted value is not already in $\supp(F)$, then the new point is that
value. If it lies in a component support $[a,b]$, then the insertion grows that
component either to $a-1$ or to $b+1$, according to the comparison between the
new decorated letter and the old component root.
\end{fact}

\begin{definition}[Witness word]
A \emph{witness word} for a forest $F$ is a decorated word $w$ such that
$\OmegaF(w)=F$.
\end{definition}

The word form of the Monk rule says that the forests appearing in
$(x_1+\cdots+x_i)P_F$ are obtained by inserting one fresh positive decorated
letter of value $i$ into a witness word for $F$ and then applying
$\OmegaF$.

\subsection{Root Lists}

When a letter is deleted, later insertions may join different pieces. Root lists let us compare the two insertion histories and track how their component boundaries change.

A root list records the components present during insertion. It starts as the
ordered list of bare number-line entries
\[
  \cdots<\underline{-1}<\underline0<\underline1<\underline2<\underline3<\cdots.
\]
As the insertion algorithm runs, each component contributes its root label to
this list, while bare entries record exposed leaves. We underline bare entries
to distinguish them from decorated letters, with
$\underline i<i^{[e]}<\underline{i+1}$ for every decoration $e$.

\begin{fact}[Root-list update]
\label{fact:root-list-update}
This update is due to Nadeau and Tewari~\cite[Lemma~5.7]{NT23}. Let $a$ be a
fresh decorated letter, and let $r_1<a<r_2$ be the two adjacent entries of the
current root list $\rl(w)$. Then
\[
  \rl(wa)=\bigl(\rl(w)\setminus\{r_1,r_2\}\bigr)\cup\{a\}.
\]
In the forest, the new vertex labeled by $a$ becomes the parent of the two
pieces represented by $r_1$ and $r_2$, where a bare integer represents a leaf.
Exactly one new support label is created by the insertion, and it is determined
by the consumed pair $\{r_1,r_2\}$.
\end{fact}

These updates let us compare the support before and after a letter is
deleted.

\subsection{Component Products and Symmetries}
\label{subsec:component-covers}

When forests have disjoint component supports, their disjoint union is written
$F\sqcup G$. On polynomials this corresponds to ordinary multiplication:
\[
  P_{F\sqcup G}=P_F P_G.
\]

The definition of an up-down forest is componentwise. We therefore
first show that the lower interval of any forest factors over its
connected components. This will reduce Theorem~B to single trees:
Boolean component intervals give a Boolean product, while a non-Boolean
interval in one component obstructs Booleanity of the whole interval.

\begin{lemma}[Component covers remain covers]
\label{lem:component-cover-remain}
Let
\[
  F=F_1\sqcup\cdots\sqcup F_r
\]
be a forest with connected components $F_j$. Suppose $E_j\lessdot F_j$ in one
component. Then
\[
  F_1\sqcup\cdots\sqcup E_j\sqcup\cdots\sqcup F_r
  \;\lessdot\;
  F_1\sqcup\cdots\sqcup F_j\sqcup\cdots\sqcup F_r .
\]
\end{lemma}

\begin{proof}
Since $E_j\lessdot F_j$, there is an index $i$ such that $P_{F_j}$ appears in
the expansion of $(x_1+\cdots+x_i)P_{E_j}$. Multiplying this expansion by the
unchanged component factors
\[
  \prod_{\ell\neq j}P_{F_\ell}
\]
shows that $P_F$ appears in the corresponding Monk product for
$F_1\sqcup\cdots\sqcup E_j\sqcup\cdots\sqcup F_r$.
\end{proof}

\begin{lemma}[Covers are component covers]
\label{lem:covers-are-component-covers}
Let
\[
  F=F_1\sqcup\cdots\sqcup F_r
\]
be a forest with connected components $F_j$, and suppose $E\lessdot F$. Let
$x$ be the unique support point in $\supp(F)\setminus\supp(E)$, and let $C_j$
be the component support of $F_j$ containing $x$. Then there is a forest $E_j$
inside $C_j\setminus\{x\}$ such that
\[
  E=F_1\sqcup\cdots\sqcup E_j\sqcup\cdots\sqcup F_r
  \qquad\text{and}\qquad
  E_j\lessdot F_j.
\]
\end{lemma}

\begin{proof}
Choose a word $w$ for $F$ whose deletion of one letter gives $E$, using
Lemma~\ref{lem:word-monk}. The gaps between the final components remain
absent in every prefix and subword, by
Theorem~\ref{thm:subword-support}, proved below. These gaps separate the
root-list histories: the letters building each final component can be
replayed independently, since no insertion crosses a gap.

Deleting one letter therefore changes only one component word; all other
components are unchanged. The two words in the affected component still
differ by one letter, so they give $E_j\lessdot F_j$. By
Fact~\ref{fact:support-growth}, this is the component containing the
unique removed support point $x$.
\end{proof}

\begin{proposition}[Component product decomposition]
\label{prop:component-product}
If $F=F_1\sqcup\cdots\sqcup F_r$ is the decomposition into connected
components, then
\[
  [\emptyforest,F]\cong
  [\emptyforest,F_1]\times\cdots\times[\emptyforest,F_r].
\]
\end{proposition}

\begin{proof}
Iterating Lemma~\ref{lem:covers-are-component-covers} down from $F$ shows that
disjoint union gives a bijection
$(E_1,\dots,E_r)\mapsto E_1\sqcup\cdots\sqcup E_r$ from the product to
$[\emptyforest,F]$. That lemma and
Lemma~\ref{lem:component-cover-remain} identify covers in both directions,
so the bijection is a poset isomorphism.
\end{proof}

Translation and reflection will allow us to treat both bent-corner
shapes through a single three-vertex example. These symmetries preserve
the Monk order. If $k\in\ZZ$, write $F+k$ for the forest obtained by
adding $k$ to every support point and every letter value, whenever all
shifted supports remain positive. If $I=[1,n]$, write $F^{\ast_I}$ for
the reflection obtained by sending $i\mapsto n+1-i$ and swapping left
and right children.

\begin{proposition}[Shift and reflection symmetries]
\label{prop:symmetries}
The following maps preserve Monk covers:
\[
  F\lessdot G\Longleftrightarrow F+k\lessdot G+k,
  \qquad
  F\lessdot G\Longleftrightarrow F^{\ast_I}\lessdot G^{\ast_I}.
\]
Consequently, whenever the shifted or reflected forests are defined,
\[
  [\emptyforest,F]\cong[\emptyforest,F+k],
  \qquad
  [\emptyforest,F]\cong[\emptyforest,F^{\ast_I}].
\]
\end{proposition}

\begin{proof}
Use witness words. Translation sends each decorated letter $i^{[e]}$ to
$(i+k)^{[e]}$. Reflection sends it to $(n+1-i)^{[-e]}$, where the negative
decoration records that equal-value tie-breakers are reversed. Translation
moves the root-list update rigidly along the number line, while reflection
reverses the root list and swaps left and right. Thus a one-letter insertion
witnessing $F\lessdot G$ is carried to a one-letter insertion witnessing the
transformed cover. Applying the inverse symmetry gives the converse.
\end{proof}

\begin{proposition}[Windows factor over runs]
\label{prop:window-factor}
Let $S\subset\ZZ_{>0}$ be finite, with maximal consecutive runs
$R_1,\dots,R_r$ of lengths $m_1,\dots,m_r$. Then
\[
  \mathcal P_S\cong
  \mathcal P_{m_1}\times\cdots\times\mathcal P_{m_r}.
\]
\end{proposition}

\begin{proof}
Restrict a forest to each run of $S$. No component crosses a gap, and
Lemmas~\ref{lem:component-cover-remain}
and~\ref{lem:covers-are-component-covers} identify covers in both
directions. Translating each run to $[1,m_j]$ by
Proposition~\ref{prop:symmetries} gives the stated product.
\end{proof}

\subsection{Words, Subwords, and Positivity}

To use a bent corner as an obstruction, we will retain selected letters
of a suitably chosen witness word. The resulting forest must still have
positive support, and we need a chain of Monk covers from it back to the
original forest. We obtain both by proving that every subword's support
is contained in the full word's support, then restoring the omitted
letters at their original positions.

We first extract the word form of the Monk rule from the insertion
correspondence of Nadeau and Tewari.

\begin{lemma}[Word form of the Monk rule]
\label{lem:word-monk}
Let $E$ be a positive-support forest and let $i\geq1$. A positive-support
forest $G$ occurs in $(x_1+\cdots+x_i)P_E$ if and only if it is obtained by
inserting a fresh decorated letter of value $i$ into some witness word for
$E$. Consequently, $E\lessdot G$ if and only if some witness word for $G$
has a one-letter deletion that is a witness for $E$.
\end{lemma}

\begin{proof}
This is the singleton-factor case of the shuffle expansion in the proof of
Nadeau--Tewari~\cite[Theorem~6.3]{NT23}, using their Proposition~5.10.
The singleton letter of value $i$ represents $x_1+\cdots+x_i$; shuffling it
with a witness class inserts it at every possible position. Positive-support
outputs give the forest summands, and the rank difference is one. In the
converse direction, a deleted letter has positive value because its value
belongs to the support of the full witness.
\end{proof}

We track a cover through the component supports of the two forests together with
the labels at their roots. For an injective word $w$, each component of
$\OmegaF(w)$ occupies a consecutive block of the number line, its \emph{component
support}; for a component support $C$ we write $\rho_w(C)$ for the decorated
label at the root of that component. The insertion map moves these blocks by the
parking rule of Nadeau and Tewari~\cite[Section~5.2]{NT23}.

\begin{lemma}[One-letter component update]
\label{lem:component-update}
Let $a$ be a fresh decorated letter, set $i=\val(a)$, and let $C=[p,q]$ be a
component support of $\OmegaF(w)$.
\begin{enumerate}[label=(\roman*)]
  \item If $i\notin\supp(\OmegaF(w))$, then the new support point is $i$. It
  joins the components, if any, ending at $i-1$ and beginning at $i+1$, and the
  new component has root label $a$.
  \item If $i\in C$, then $C$ is enlarged to the left when $a<\rho_w(C)$ and to
  the right when $a>\rho_w(C)$. The new support point is respectively $p-1$ or
  $q+1$, and any component adjacent to that new point is joined as well. The
  resulting component has root label $a$.
\end{enumerate}
Moreover, root labels increase from left to right among the component supports.
\end{lemma}

\begin{proof}
The support update is the parking procedure of Nadeau and
Tewari~\cite[Section~5.2]{NT23}, and the description of the new root is the
insertion definition of $\OmegaF$. A root label has value in its component
support, by induction on the insertions, so the root labels of disjoint
component supports occur in their left-to-right order.
\end{proof}

Support containment alone does not record the root comparisons that
determine the direction of the next insertion. To prove that containment
persists along a word, we strengthen it by comparing component roots at
shared endpoints.

\begin{definition}[Boundary-root refinement]
\label{def:boundary-refinement}
For two labeled forests $A$ and $B$, write $A\preccurlyeq B$ when the following
hold.
\begin{enumerate}[label=(\alph*)]
  \item Every component support $C$ of $A$ is contained in a component support
  $D$ of $B$.
  \item If $C$ and $D$ have the same left endpoint, then
  $\rho_A(C)\leq\rho_B(D)$.
  \item If $C$ and $D$ have the same right endpoint, then
  $\rho_A(C)\geq\rho_B(D)$.
\end{enumerate}
\end{definition}

The endpoint inequalities force the two parking procedures to agree on their
direction of growth. If a component of $A$ reaches the right boundary of its
containing component of $B$, then a letter larger than the smaller root is also
larger than the larger root, so both procedures enlarge to the right whenever
the smaller one would leave the containing support; the left case is symmetric.

\begin{lemma}[Stability under one insertion]
\label{lem:insertion-stability}
Suppose $A\preccurlyeq B$, and let $a$ be a fresh letter. Then
\[
  A\preccurlyeq Ba
  \qquad\text{and}\qquad
  Aa\preccurlyeq Ba,
\]
where $Ca$ denotes the result of inserting $a$ into the labeled forest $C$.
\end{lemma}

\begin{proof}
The boundary inequalities prevent a smaller component from growing past
its containing component. Set $i=\val(a)$ and use
Lemma~\ref{lem:component-update}; only components containing or neighboring
$i$ can change.

First suppose $i\notin\supp(B)$. Then also $i\notin\supp(A)$. Both updates
create the support point $i$ and join only components immediately adjacent to
it. Every newly joined component of $A$ is contained in the corresponding newly
joined component of $B$. At any outer endpoint inherited from a left component,
the old inequality is strengthened because the new root is $a$ and the root of
that left component is less than $a$. The right-endpoint statement is the
reflected argument. This proves both refinements in this case.

\begin{figure}[H]
\centering
\begin{tikzpicture}[
  x=.72cm,y=.72cm,
  line cap=round,line join=round,
  forestedge/.style={line width=.85pt},
  newedge/.style={line width=1pt,green!55!black},
  dot/.style={circle,fill=black,inner sep=1.45pt},
  newdot/.style={circle,fill=green!55!black,inner sep=1.55pt},
  lab/.style={red!80!black,font=\scriptsize,fill=white,inner sep=.6pt},
  newlab/.style={green!55!black,font=\scriptsize,fill=white,inner sep=.6pt}
]
  \begin{scope}[shift={(0,0)}]
    \forestline{1}{6}
    \coordinate (l1) at (1.5,.62);
    \coordinate (l) at (2.25,1.20);
    \coordinate (r1) at (4.5,.62);
    \coordinate (r) at (5.25,1.20);
    \draw[forestedge] (l1)--(1,0) (l1)--(2,0);
    \draw[forestedge] (l)--(l1) (l)--(3,0);
    \draw[forestedge] (r1)--(4,0) (r1)--(5,0);
    \draw[forestedge] (r)--(r1) (r)--(6,0);
    \node[dot] at (l1) {}; \node[dot] at (l) {};
    \node[dot] at (r1) {}; \node[dot] at (r) {};
    \node[lab] at (2.55,1.45) {$r_L$};
    \node[lab] at (5.55,1.45) {$r_R$};
  \end{scope}

  \draw[-{Stealth[length=2mm]},line width=.8pt] (7.2,1.0)--(9.1,1.0);

  \begin{scope}[shift={(10,0)}]
    \forestline{1}{6}
    \coordinate (l1) at (1.5,.62);
    \coordinate (l) at (2.25,1.20);
    \coordinate (r1) at (4.5,.62);
    \coordinate (r) at (5.25,1.20);
    \coordinate (a) at (3.75,2.05);
    \draw[forestedge] (l1)--(1,0) (l1)--(2,0);
    \draw[forestedge] (l)--(l1) (l)--(3,0);
    \draw[forestedge] (r1)--(4,0) (r1)--(5,0);
    \draw[forestedge] (r)--(r1) (r)--(6,0);
    \draw[newedge] (a)--(l) (a)--(r);
    \node[dot] at (l1) {}; \node[dot] at (l) {};
    \node[dot] at (r1) {}; \node[dot] at (r) {};
    \node[newdot] at (a) {};
    \node[lab] at (2.55,1.45) {$r_L$};
    \node[lab] at (5.55,1.45) {$r_R$};
    \node[newlab] at (3.75,2.35) {$a$};
  \end{scope}
\end{tikzpicture}
\caption{$i$ not in support}
\label{fig:subword-case-empty}
\end{figure}
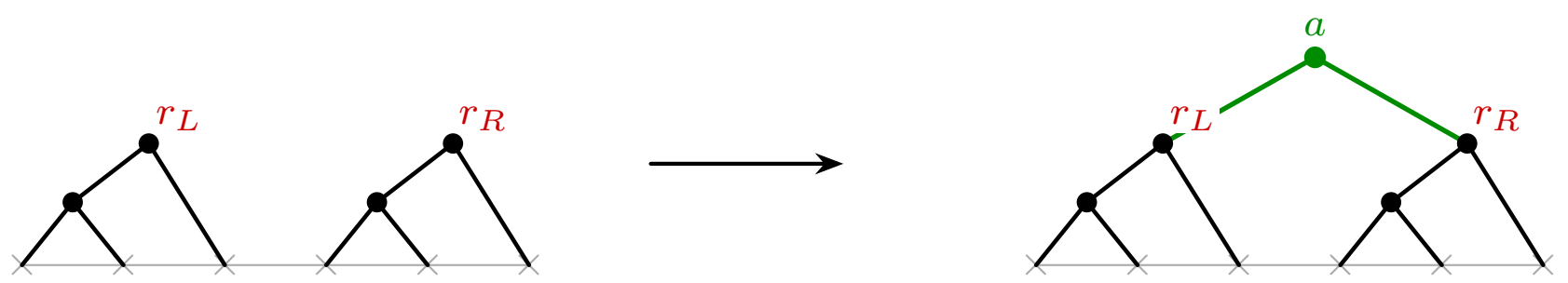

If $i\in\supp(B)\setminus\supp(A)$, the new $A$-component lies inside
the old $B$-component containing $i$, while the $B$-update can only
enlarge that component. At a shared endpoint, an unchanged $A$-root
retains its inequality by $A\preccurlyeq B$ and the root order; a newly
created $A$-root is $a$, the same as the new $B$-root. Thus both
refinements hold.

\begin{figure}[H]
\centering
\begin{tikzpicture}[
  x=.62cm,y=.62cm,
  line cap=round,line join=round,
  forestedge/.style={line width=.85pt},
  newedge/.style={line width=1pt,green!55!black},
  dot/.style={circle,fill=black,inner sep=1.35pt},
  newdot/.style={circle,fill=green!55!black,inner sep=1.5pt},
  lab/.style={red!80!black,font=\scriptsize,fill=white,inner sep=.55pt},
  newlab/.style={green!55!black,font=\scriptsize,fill=white,inner sep=.55pt}
]
  \node[anchor=east,font=\small] at (.55,2.85) {$B$};
  \node[anchor=east,font=\small] at (.55,-2.75) {$A$};

  \begin{scope}[shift={(0,2.2)}]
    \forestline{1}{5}
    \coordinate (b1) at (1.5,.60);
    \coordinate (b2) at (2.25,1.18);
    \coordinate (broot) at (3.12,1.72);
    \draw[forestedge] (b1)--(1,0) (b1)--(2,0);
    \draw[forestedge] (b2)--(b1) (b2)--(3,0);
    \draw[forestedge] (broot)--(b2) (broot)--(4,0);
    \node[dot] at (b1) {}; \node[dot] at (b2) {};
    \node[dot] at (broot) {};
    \node[lab] at (3.45,2.00) {$r_B$};
  \end{scope}

  \draw[-{Stealth[length=2mm]},line width=.8pt] (6.0,3.03)--(8.05,3.03);

  \begin{scope}[shift={(10,2.2)}]
    \forestline{1}{5}
    \coordinate (b1) at (1.5,.60);
    \coordinate (b2) at (2.25,1.18);
    \coordinate (broot) at (3.12,1.72);
    \coordinate (a) at (4.12,2.25);
    \draw[forestedge] (b1)--(1,0) (b1)--(2,0);
    \draw[forestedge] (b2)--(b1) (b2)--(3,0);
    \draw[forestedge] (broot)--(b2) (broot)--(4,0);
    \draw[newedge] (a)--(broot) (a)--(5,0);
    \node[dot] at (b1) {}; \node[dot] at (b2) {};
    \node[dot] at (broot) {};
    \node[newdot] at (a) {};
    \node[lab] at (3.45,2.00) {$r_B$};
    \node[newlab] at (4.12,2.55) {$a$};
  \end{scope}

  \begin{scope}[shift={(0,-3.4)}]
    \forestline{1}{5}
    \coordinate (al) at (1.5,.66);
    \coordinate (ar) at (3.5,.66);
    \draw[forestedge] (al)--(1,0) (al)--(2,0);
    \draw[forestedge] (ar)--(3,0) (ar)--(4,0);
    \node[dot] at (al) {};
    \node[dot] at (ar) {};
  \end{scope}

  \draw[-{Stealth[length=2mm]},line width=.8pt] (6.0,-2.57)--(8.05,-2.57);

  \begin{scope}[shift={(10,-3.4)}]
    \forestline{1}{5}
    \coordinate (al) at (1.5,.66);
    \coordinate (ar) at (3.5,.66);
    \coordinate (a) at (2.5,1.52);
    \draw[forestedge] (al)--(1,0) (al)--(2,0);
    \draw[forestedge] (ar)--(3,0) (ar)--(4,0);
    \draw[newedge] (a)--(al) (a)--(ar);
    \node[dot] at (al) {};
    \node[dot] at (ar) {};
    \node[newdot] at (a) {};
    \node[newlab] at (2.5,1.88) {$a$};
  \end{scope}
\end{tikzpicture}
\caption{$i\in\supp(B)\setminus\supp(A)$.}
\label{fig:subword-case-middle}
\end{figure}
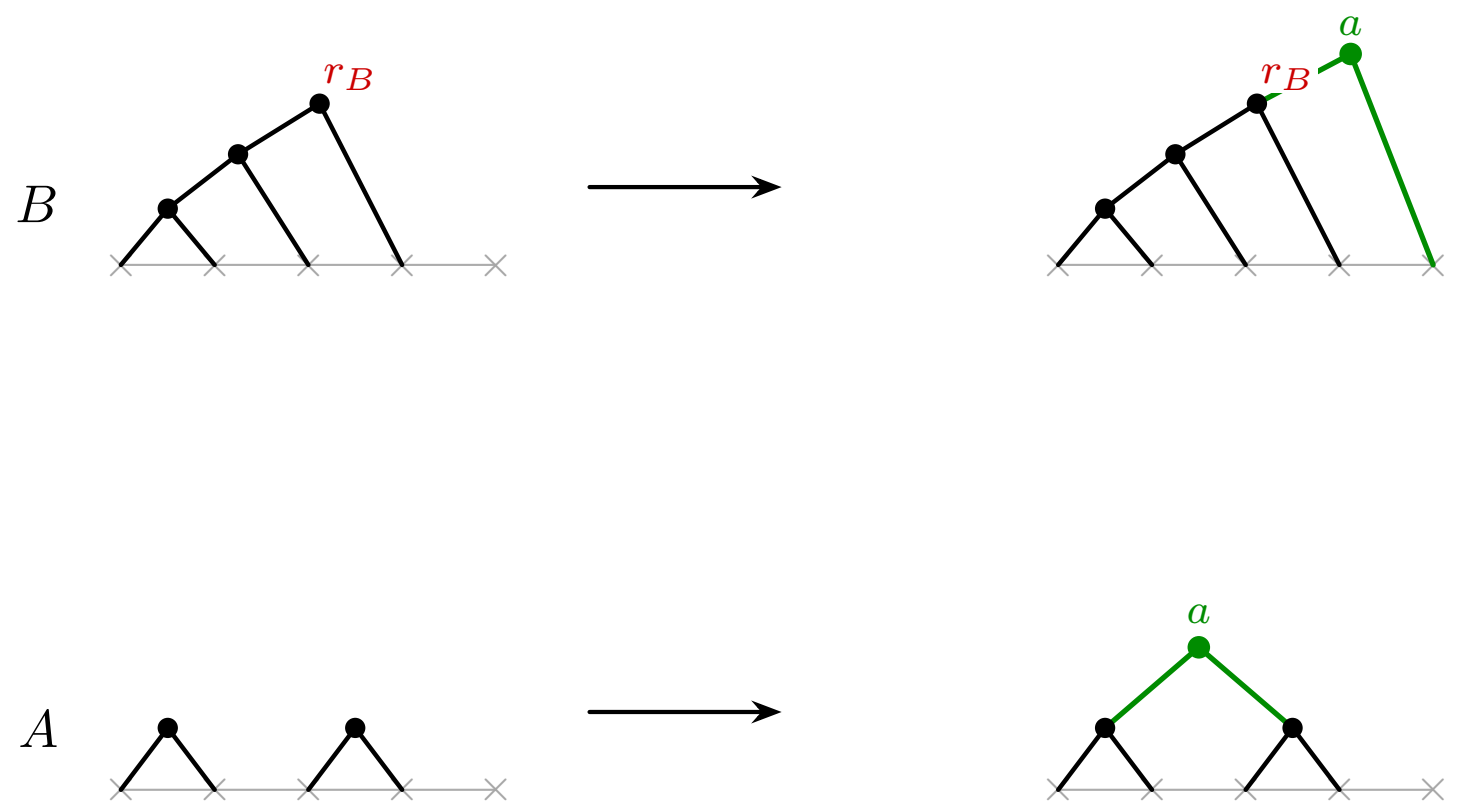

Finally, suppose $i\in\supp(A)$, and let $C\subseteq D$ be the component
supports containing $i$. Growth of $C$ can leave $D$ only through a
shared endpoint. On the right, this requires
\[
  a>\rho_A(C)\geq\rho_B(D),
\]
so the $B$-component grows right as well; the left case is its reflection.
This proves containment. Both new roots are $a$, unchanged outer endpoints
retain their inequalities, and any absorbed neighboring component uses
the left-to-right root order. Hence $Aa\preccurlyeq Ba$. The same
endpoint comparisons with $A$ left unchanged give $A\preccurlyeq Ba$. \qedhere

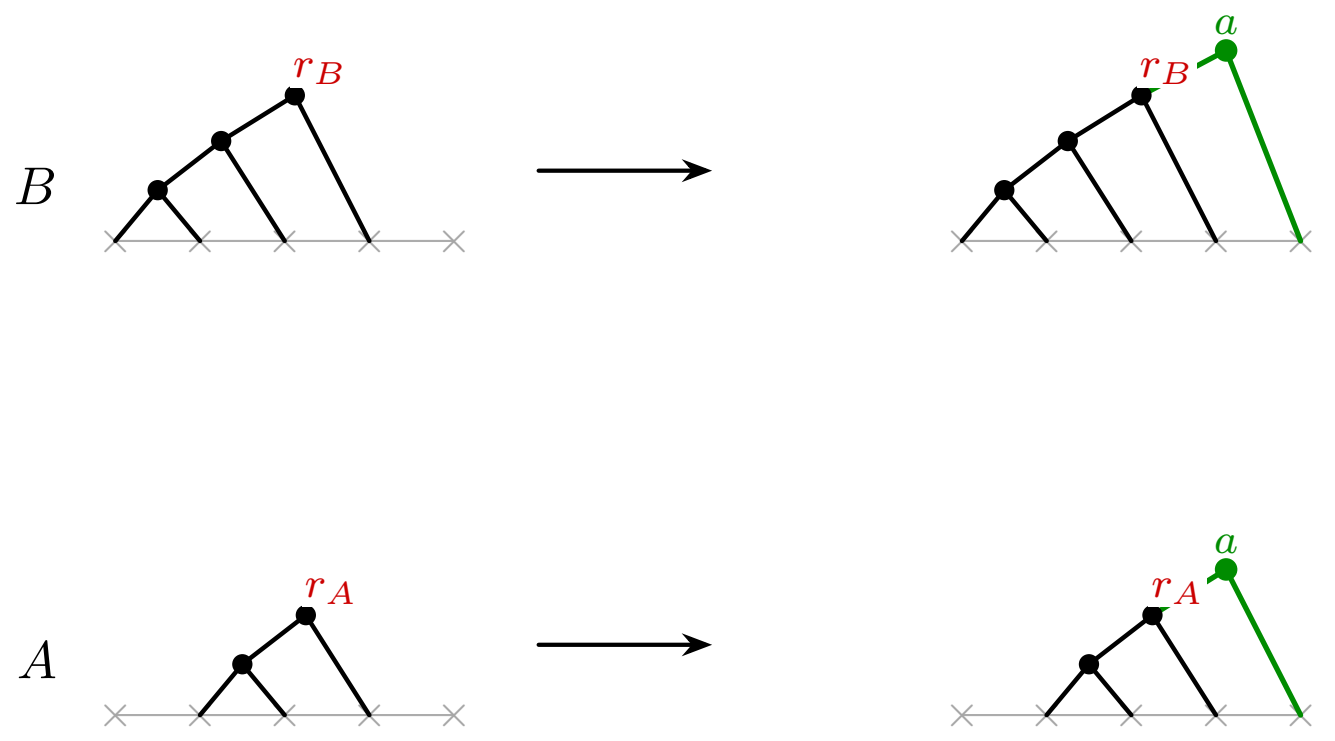
\begin{figure}[H]
\centering
\begin{tikzpicture}[
  x=.56cm,y=.56cm,
  line cap=round,line join=round,
  forestedge/.style={line width=.85pt},
  newedge/.style={line width=1pt,green!55!black},
  dot/.style={circle,fill=black,inner sep=1.35pt},
  newdot/.style={circle,fill=green!55!black,inner sep=1.5pt},
  lab/.style={red!80!black,font=\scriptsize,fill=white,inner sep=.55pt},
  newlab/.style={green!55!black,font=\scriptsize,fill=white,inner sep=.55pt}
]
  \node[anchor=east,font=\small] at (.55,2.85) {$B$};
  \node[anchor=east,font=\small] at (.55,-2.75) {$A$};

  \begin{scope}[shift={(0,2.2)}]
    \forestline{1}{5}
    \coordinate (b1) at (1.5,.60);
    \coordinate (b2) at (2.25,1.18);
    \coordinate (broot) at (3.12,1.72);
    \draw[forestedge] (b1)--(1,0) (b1)--(2,0);
    \draw[forestedge] (b2)--(b1) (b2)--(3,0);
    \draw[forestedge] (broot)--(b2) (broot)--(4,0);
    \node[dot] at (b1) {}; \node[dot] at (b2) {};
    \node[dot] at (broot) {};
    \node[lab] at (3.42,2.00) {$r_B$};
  \end{scope}

  \draw[-{Stealth[length=2mm]},line width=.8pt] (6.0,3.03)--(8.05,3.03);

  \begin{scope}[shift={(10,2.2)}]
    \forestline{1}{5}
    \coordinate (b1) at (1.5,.60);
    \coordinate (b2) at (2.25,1.18);
    \coordinate (broot) at (3.12,1.72);
    \coordinate (a) at (4.12,2.25);
    \draw[forestedge] (b1)--(1,0) (b1)--(2,0);
    \draw[forestedge] (b2)--(b1) (b2)--(3,0);
    \draw[forestedge] (broot)--(b2) (broot)--(4,0);
    \draw[newedge] (a)--(broot) (a)--(5,0);
    \node[dot] at (b1) {}; \node[dot] at (b2) {};
    \node[dot] at (broot) {};
    \node[newdot] at (a) {};
    \node[lab] at (3.42,2.00) {$r_B$};
    \node[newlab] at (4.12,2.55) {$a$};
  \end{scope}

  \begin{scope}[shift={(0,-3.4)}]
    \forestline{1}{5}
    \coordinate (a1) at (2.5,.60);
    \coordinate (aroot) at (3.25,1.18);
    \draw[forestedge] (a1)--(2,0) (a1)--(3,0);
    \draw[forestedge] (aroot)--(a1) (aroot)--(4,0);
    \node[dot] at (a1) {}; \node[dot] at (aroot) {};
    \node[lab] at (3.55,1.46) {$r_A$};
  \end{scope}

  \draw[-{Stealth[length=2mm]},line width=.8pt] (6.0,-2.57)--(8.05,-2.57);

  \begin{scope}[shift={(10,-3.4)}]
    \forestline{1}{5}
    \coordinate (a1) at (2.5,.60);
    \coordinate (aroot) at (3.25,1.18);
    \coordinate (a) at (4.12,1.72);
    \draw[forestedge] (a1)--(2,0) (a1)--(3,0);
    \draw[forestedge] (aroot)--(a1) (aroot)--(4,0);
    \draw[newedge] (a)--(aroot) (a)--(5,0);
    \node[dot] at (a1) {}; \node[dot] at (aroot) {};
    \node[newdot] at (a) {};
    \node[lab] at (3.55,1.46) {$r_A$};
    \node[newlab] at (4.12,2.02) {$a$};
  \end{scope}
\end{tikzpicture}
\caption{The shared-right-boundary case.}
\label{fig:subword-case-boundary}
\end{figure}
\end{proof}

The refinement propagates along a word one letter at a time, so it holds for any
subword.

\begin{proposition}[Subword refinement]
\label{prop:subword-refinement}
Let $u$ be a subword of $w$. Then $\OmegaF(u)\preccurlyeq\OmegaF(w)$.
\end{proposition}

\begin{proof}
Read $w$ from left to right, starting with two empty forests. At each letter,
apply the first assertion of Lemma~\ref{lem:insertion-stability} if the letter
is omitted from $u$, and the second if it is retained. Induction gives
$\OmegaF(u)\preccurlyeq\OmegaF(w)$.
\end{proof}

\begin{theorem}[Subword support containment]
\label{thm:subword-support}
Let $u$ be a subword of an injective word $w$ in decorated letters. Then
\[
  \supp(\OmegaF(u))\subseteq\supp(\OmegaF(w)).
\]
Consequently, if $w$ has positive support, then every subword of $w$ has
positive support as well.
\end{theorem}

\begin{proof}
By Proposition~\ref{prop:subword-refinement} and
Definition~\ref{def:boundary-refinement}(a), every component support of
$\OmegaF(u)$ lies in a component support of $\OmegaF(w)$. Taking the union over
components gives the containment $\supp(\OmegaF(u))\subseteq\supp(\OmegaF(w))$.
If the larger support is contained in $\ZZ_{>0}$, then so is the smaller.
\end{proof}

\begin{corollary}[Arbitrary restoration order]
\label{cor:restoration}
Let $u$ be a subword of a positive-support witness word $w$. Restore the letters
of $w\setminus u$ in any order, always at their original positions. Every
intermediate word has positive support, and every restoration is a valid
one-letter Monk cover.
\end{corollary}

\begin{proof}
Each intermediate word is a subword of $w$, so it has positive support by
Theorem~\ref{thm:subword-support}. Restoring one letter is a positive-valued
one-letter insertion into a witness word, so the word form of the Monk rule
(Lemma~\ref{lem:word-monk}) gives the required cover.
\end{proof}

The restoration corollary will be used in
Theorem~\ref{thm:zigzag-isolation} to place a bent-corner shape below the
original forest. We now establish another ingredient for Theorem~B:
every prescribed support point can be removed by a lower cover.

\subsection{Support-Removing Covers}

Removing points successively will produce a lower element for every
subset of $\supp(F)$. We first prove the existence of a cover removing
any one point, using a witness whose letter values are precisely the
support positions. For an up-down forest, Section~\ref{sec:boolean}
will add the crucial uniqueness statement: the remaining support
determines the lower cover.

\begin{theorem}[Support-removing covers]
\label{thm:support-removing}
Let $F$ be an indexed forest and let $x\in\supp(F)$. Then there is a forest
$F^{\setminus x}$ such that
\[
  F^{\setminus x}\lessdot F
  \qquad\text{and}\qquad
  \supp(F^{\setminus x})=\supp(F)\setminus\{x\}.
\]
\end{theorem}

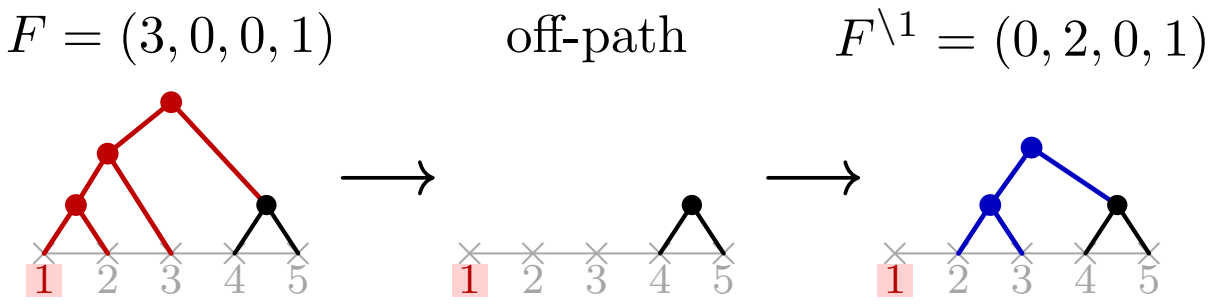
\begin{figure}[H]
\centering
\begin{tikzpicture}[
  x=.42cm, y=.40cm,
  line cap=round, line join=round,
  forestedge/.style={line width=.8pt},
  rededge/.style={line width=.9pt,red!75!black},
  blueedge/.style={line width=.9pt,blue!75!black},
  dot/.style={circle,fill=black,inner sep=1.35pt},
  reddot/.style={circle,fill=red!75!black,inner sep=1.45pt},
  bluedot/.style={circle,fill=blue!75!black,inner sep=1.45pt},
  leaf/.style={gray!70,font=\scriptsize},
  title/.style={font=\small},
  arr/.style={->,line width=.7pt}]
  \newcommand{\srline}{%
    \draw[gray!65] (0,0)--(4,0);
    \foreach \x/\lab in {0/1,1/2,2/3,3/4,4/5} {
      \node[leaf] at (\x,0) {$\times$};
      \node[leaf] at (\x,-.42) {\lab};
    }
    \fill[red!18] (-.28,-.72) rectangle (.28,-.18);
    \node[font=\scriptsize,red!70!black] at (0,-.42) {$1$};
  }

  \begin{scope}[shift={(0,0)}]
    \node[title] at (2,3.55) {$F=(3,0,0,1)$};
    \srline
    \draw[forestedge] (3.5,.8)--(3,0) (3.5,.8)--(4,0);
    \draw[rededge]
      (2.0,2.5)--(1.0,1.65) (2.0,2.5)--(3.5,.8)
      (1.0,1.65)--(.5,.8) (1.0,1.65)--(2,0)
      (.5,.8)--(0,0) (.5,.8)--(1,0);
    \node[reddot] at (.5,.8) {};
    \node[reddot] at (1.0,1.65) {};
    \node[reddot] at (2.0,2.5) {};
    \node[dot] at (3.5,.8) {};
  \end{scope}

  \draw[arr] (4.7,1.25)--(6.15,1.25);

  \begin{scope}[shift={(6.7,0)}]
    \node[title] at (2,3.55) {off-path};
    \srline
    \draw[forestedge] (3.5,.8)--(3,0) (3.5,.8)--(4,0);
    \node[dot] at (3.5,.8) {};
  \end{scope}

  \draw[arr] (11.4,1.25)--(12.85,1.25);

  \begin{scope}[shift={(13.4,0)}]
    \node[title] at (2,3.55) {$F^{\setminus 1}=(0,2,0,1)$};
    \srline
    \draw[forestedge] (3.5,.8)--(3,0) (3.5,.8)--(4,0);
    \draw[blueedge]
      (2.15,1.75)--(1.5,.8) (2.15,1.75)--(3.5,.8)
      (1.5,.8)--(1,0) (1.5,.8)--(2,0);
    \node[bluedot] at (1.5,.8) {};
    \node[bluedot] at (2.15,1.75) {};
    \node[dot] at (3.5,.8) {};
  \end{scope}
\end{tikzpicture}
\caption{Removing the support point $x=1$ from
$(3,0,0,1)$.}
\label{fig:support-removing}
\end{figure}

For the tree in Figure~\ref{fig:support-removing}, deleting $1^{[1]}$ from
the witness $1^{[1]}2^{[1]}4^{[1]}3^{[1]}$ gives the displayed lower cover.
The middle panel shows the cherry unchanged by this deletion.

\begin{proof}[Proof of Theorem~\ref{thm:support-removing}]
Suppose first that $F$ is connected, and label its internal vertices by
their support positions in inorder. For a subtree $U$ with root label $r$,
define the postorder word
\[
  w(U)=w(U_L)\,w(U_R)\,r^{[1]},
\]
with empty subtrees contributing empty words. Inductively,
$\OmegaF(w(U))=U$: the first two factors build the subtrees on opposite
sides of the absent position $r$, and inserting $r^{[1]}$ joins them
under their prescribed parent, using a bare leaf for an empty subtree.

All values in $w(F)$ are distinct. Consequently, for every subword $u$,
\[
  \supp(\OmegaF(u))=\{i:i^{[1]}\text{ occurs in }u\}.
\]
Indeed, each successive value is absent from the current support, so its
insertion adjoins exactly that value. Delete $x^{[1]}$ from $w(F)$ to
obtain $h$. The displayed identity gives
$\supp(\OmegaF(h))=\supp(F)\setminus\{x\}$, in particular positive support.
Restoring $x^{[1]}$ at its original position recovers $w(F)$ exactly, so
Lemma~\ref{lem:word-monk} gives $\OmegaF(h)\lessdot F$.

For disconnected $F$, apply the same construction in the component
containing $x$ and leave the others unchanged.
Lemma~\ref{lem:component-cover-remain} extends its cover to $F$.
\end{proof}

\subsection{Forest Division}

The insertion model also explains where a lower cover can change a
forest. Deleting a letter alters attachments only along its ancestor
path; the subtrees off that path remain unchanged. The final subsections
factor out this common base to describe the cover through a path-shaped
quotient.

The forest product and quotient are defined by
Nadeau--Spink--Tewari~\cite[Section~4.1]{NST26}. If $Q$ and $B$ are forests,
their product $Q\cdot B$ is obtained by grafting the $i$th root of $B$ onto
the $i$th leaf of $Q$. Equivalently, factorizations
\[
  H=Q\cdot B
\]
are the same as decompositions of the internal vertices of $H$ into two
parts: the quotient part $Q$, which is closed upward under taking parents, and
the base part $B$, which is closed downward under taking children. The product
is right-cancellative~\cite[Proposition~4.3]{NST26}, so if $H=Q\cdot B$ then $Q$ is unique; we write
\[
  Q=H/B.
\]
The common-base construction below uses this quotient.
\subsection{Common Bases and Spine Locality}

Deleting a letter need not simply erase one vertex from the final drawing:
later insertions can attach differently. The changes stay above that letter,
while the subtrees beside its ancestor path remain intact.

\begin{definition}[Affected spine and base]
\label{def:affected-spine}
Let $w=w_-aw_+$ be an injective decorated word, and put
$\widehat w=w_-w_+$. In the labeled forest for $w$, let $\Spine_a(w)$ be
the parent chain from $a$ to the root of its component, including $a$.
The complementary vertices form the \emph{base} $\Base_a(w)$. This choice
depends on the witness and may give the empty base.
\end{definition}

\begin{lemma}[Spine locality]
\label{lem:spine-locality}
Let $w=w_-aw_+$ and $\widehat w=w_-w_+$ as above, and suppose both
$\OmegaF(w)$ and $\OmegaF(\widehat w)$ have positive support. Every
vertex in $\Base_a(w)$ has the same two children in the two labeled
forests.
\end{lemma}

\begin{proof}
Run the two root-list histories together. Immediately after $a$ is
inserted, they agree except that the upper history has the root $a$ where
the lower history has the two pieces consumed by $a$. The same description
persists through the common suffix: one upper root $r$ is replaced by two
consecutive lower entries $u,v$ with $u\leq r\leq v$.

A later insertion avoiding $r$ consumes the same two entries in both
histories. An insertion consuming $r$ becomes its parent in the upper
history; replacing the corresponding adjacent pair in the lower history
again leaves two entries bracketing the new upper root. Thus only the
successive ancestors of $a$ can acquire different children. All vertices
off that path retain their children and the subtrees below them.
\end{proof}

\begin{lemma}[Lower covers have a common base]
\label{lem:common-base}
If $F\lessdot G$ is realized by inserting the letter $a$, then $F$ and
$G$ share a child-closed base forest $B$ whose complement in $G$ is the
affected spine above $a$.
\end{lemma}

\begin{proof}
Choose witnesses $w=w_-aw_+$ for $G$ and $\widehat w=w_-w_+$ for $F$,
and take the vertices outside $\Spine_a(w)$. The complement of an
ancestor path is child-closed, and Lemma~\ref{lem:spine-locality} shows
that these subtrees occur unchanged in both forests. They form $B$.
\end{proof}

For example, deleting the first letter of $1^{[1]}1^{[2]}$ gives the
cover $(1)\lessdot(2)$. Both upper vertices lie on the affected spine,
so the chosen base is empty and the quotient is the cover itself.

\subsection{Spine Quotients}

Contracting the unchanged base leaves a cover of a path-shaped tree.
The contraction is expressed by forest division.

\begin{definition}[Spine quotient by forest division]
\label{def:spine-quotient}
Let $F\lessdot G$ be realized by $w=w_-aw_+$ and $\widehat w=w_-w_+$,
and let $B$ be the base from Lemma~\ref{lem:common-base}. Since $B$ is
child-closed in both forests, write
\[
  G=Z\cdot B,\qquad F=H\cdot B.
\]
The pair $H=F/B$, $Z=G/B$ is the \emph{spine quotient} of this witness.
\end{definition}

\begin{lemma}[Quotient covers remain covers]
\label{lem:quotient-cover}
With notation as in Definition~\ref{def:spine-quotient},
\[
  H\lessdot Z.
\]
Moreover, $Z$ is a possibly shifted path-shaped tree with the same left/right turn word
as the affected spine of $G$.
\end{lemma}

\begin{proof}
Let $b$ be the base letters and $z$ the remaining letters of $w$, each in
its original order. Since the base is child-closed, we may build it first:
every remaining vertex still finds its prescribed children when reached.
Thus $bz$ and $b\widehat z$ are witnesses for $G$ and $F$, where
$\widehat z$ is obtained from $z$ by deleting $a$.

After building $B$, contract its ordered root-list entries to bare quotient
leaves. Relabel the letters of $z$ between these entries, preserving their
order. Root-list insertion depends only on order and adjacency, so the
two resulting histories construct exactly $Z$ and $H$. Their words still
differ by the single letter corresponding to $a$; hence
Lemma~\ref{lem:word-monk} gives $H\lessdot Z$.
The upper quotient retains precisely the ancestor path, with its original
left and right attachments.
\end{proof}

\begin{theorem}[All lower covers reduce to zigzag quotients]
\label{thm:all-covers-zigzag}
Let $F\lessdot G$ be a lower cover. There is a common child-closed base
$B$ such that
\[
  F=H\cdot B,\qquad G=Z\cdot B,\qquad H\lessdot Z,
\]
where $Z$ is a possibly shifted path-shaped tree. Conversely, grafting the same base onto
a cover of such a quotient gives a lower cover.
\end{theorem}

\begin{proof}
The forward direction is Lemmas~\ref{lem:common-base}
and~\ref{lem:quotient-cover}. For the converse, choose one-letter-related
witnesses for $H\lessdot Z$ and build $B$ first, choosing its decorations
with enough room for the finitely many lifted letters. Replace each bare
quotient leaf by the corresponding root-list entry of $B$, lifting the
remaining letters in the same order between consecutive entries. This
reverses the contraction above, so the histories construct $H\cdot B$
and $Z\cdot B$ and still differ by one letter. The word form of the Monk
rule gives the grafted cover.
\end{proof}

We now apply the insertion tools to Theorem~B. For up-down forests,
support-removing covers will give every support subset, and an analysis
of witness-word deletions will establish uniqueness. For all other
forests, subword restoration will place a bent-corner obstruction below
the top.

\section{Boolean Lower Intervals}
\label{sec:boolean}

The principal lower intervals that are Boolean lattices are exactly those
below disjoint unions of up-down trees.

\begin{definition}[Up-down trees]
  \label{def:updown-tree}
  An \emph{up-down tree} is a tree whose internal vertices consist of the
  root, a straight branch of left children, and a straight branch of right
  children; either branch may have length zero. An \emph{up-down forest} is
  a disjoint union of up-down trees.
\end{definition}

\begin{theorem}[Boolean interval classification]
\label{thm:boolean-classification}
For a forest $F$, the lower interval $[\emptyforest,F]$ is Boolean if and only
if $F$ is an up-down forest. In that case
\[
  [\emptyforest,F]\cong B_{|\supp(F)|}.
\]
\end{theorem}

In code coordinates, the up-down tree with support $[a,b]$ and turn $t$, where
$a\leq t\leq b$, has
\[
  c_a=t-a+1,\qquad
  c_{a+1}=\cdots=c_t=0,\qquad
  c_{t+1}=\cdots=c_b=1,
\]
with all other code entries equal to zero. Thus $(2,0,1,1)$ is the up-down tree
on $[1,4]$ with turn $2$, while $(3,0,0,1)$ is the up-down tree on $[1,4]$
with turn $3$.

\subsection{A First Boolean Interval}

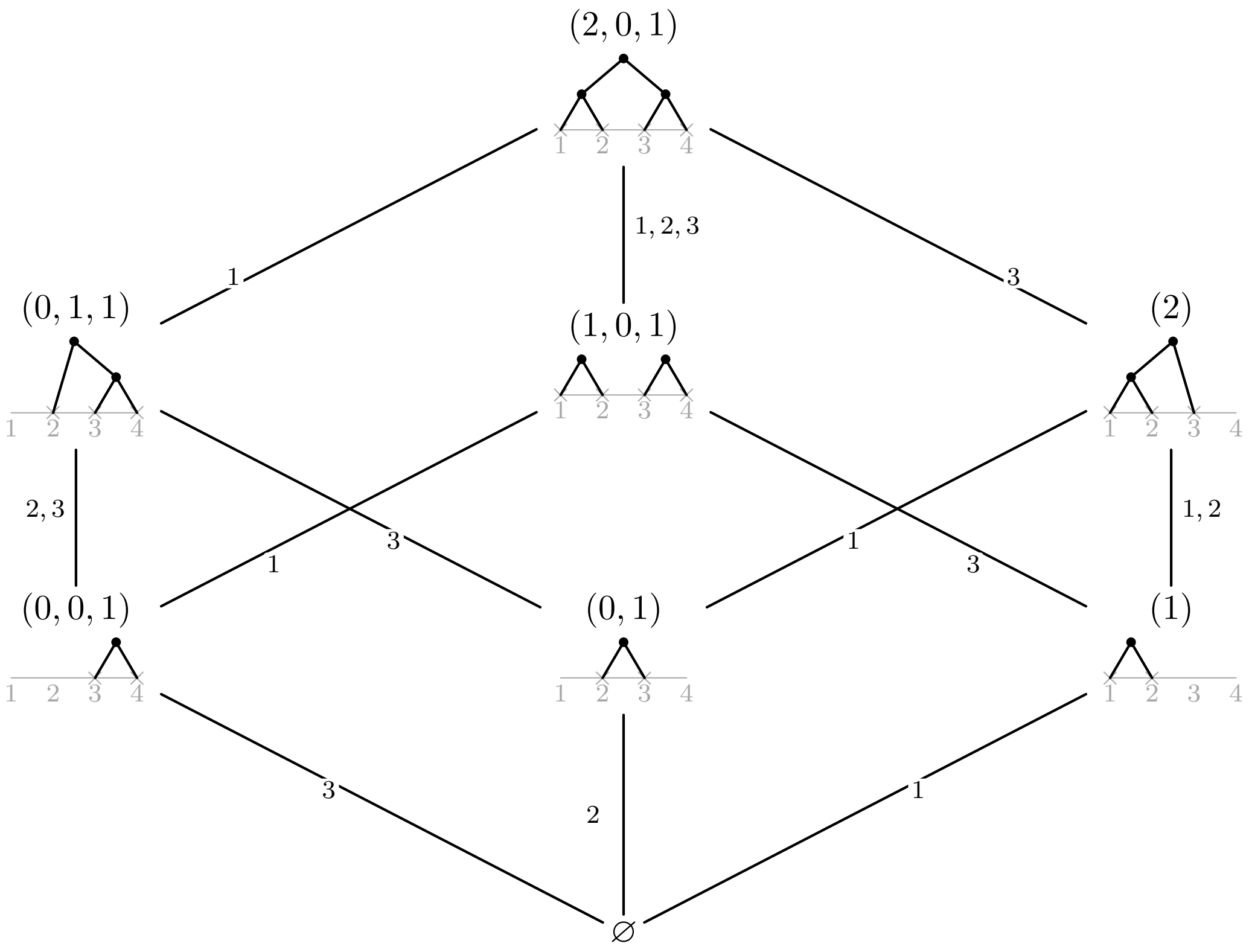
\begin{figure}[H]
\centering
\begingroup

\newcommand{\fOne}{%
\begin{tikzpicture}[x=.46cm,y=.46cm,baseline=-0.6ex,line cap=round,line join=round]
    \draw[gray!65] (0,0)--(3,0);
    \foreach \x/\lab in {0/1,1/2,2/3,3/4} {
        \node[below,gray!70] at (\x,0) {\scriptsize \lab};
    }
    \node[gray!70] at (0,0) {\scriptsize $\times$};
    \node[gray!70] at (1,0) {\scriptsize $\times$};
    \draw[thick] (0.5,.85)--(0,0);
    \draw[thick] (0.5,.85)--(1,0);
    \fill (0.5,.85) circle (1.5pt);
\end{tikzpicture}}

\newcommand{\fTwo}{%
\begin{tikzpicture}[x=.46cm,y=.46cm,baseline=-0.6ex,line cap=round,line join=round]
    \draw[gray!65] (0,0)--(3,0);
    \foreach \x/\lab in {0/1,1/2,2/3,3/4} {
        \node[below,gray!70] at (\x,0) {\scriptsize \lab};
    }
    \node[gray!70] at (1,0) {\scriptsize $\times$};
    \node[gray!70] at (2,0) {\scriptsize $\times$};
    \draw[thick] (1.5,.85)--(1,0);
    \draw[thick] (1.5,.85)--(2,0);
    \fill (1.5,.85) circle (1.5pt);
\end{tikzpicture}}

\newcommand{\fThree}{%
\begin{tikzpicture}[x=.46cm,y=.46cm,baseline=-0.6ex,line cap=round,line join=round]
    \draw[gray!65] (0,0)--(3,0);
    \foreach \x/\lab in {0/1,1/2,2/3,3/4} {
        \node[below,gray!70] at (\x,0) {\scriptsize \lab};
    }
    \node[gray!70] at (2,0) {\scriptsize $\times$};
    \node[gray!70] at (3,0) {\scriptsize $\times$};
    \draw[thick] (2.5,.85)--(2,0);
    \draw[thick] (2.5,.85)--(3,0);
    \fill (2.5,.85) circle (1.5pt);
\end{tikzpicture}}

\newcommand{\fOneThree}{%
\begin{tikzpicture}[x=.46cm,y=.46cm,baseline=-0.6ex,line cap=round,line join=round]
    \draw[gray!65] (0,0)--(3,0);
    \foreach \x/\lab in {0/1,1/2,2/3,3/4} {
        \node[below,gray!70] at (\x,0) {\scriptsize \lab};
    }
    \node[gray!70] at (0,0) {\scriptsize $\times$};
    \node[gray!70] at (1,0) {\scriptsize $\times$};
    \node[gray!70] at (2,0) {\scriptsize $\times$};
    \node[gray!70] at (3,0) {\scriptsize $\times$};
    \draw[thick] (0.5,.85)--(0,0);
    \draw[thick] (0.5,.85)--(1,0);
    \fill (0.5,.85) circle (1.5pt);
    \draw[thick] (2.5,.85)--(2,0);
    \draw[thick] (2.5,.85)--(3,0);
    \fill (2.5,.85) circle (1.5pt);
\end{tikzpicture}}

\newcommand{\fZeroOneOne}{%
\begin{tikzpicture}[x=.46cm,y=.46cm,baseline=-0.6ex,line cap=round,line join=round]
    \draw[gray!65] (0,0)--(3,0);
    \foreach \x/\lab in {0/1,1/2,2/3,3/4} {
        \node[below,gray!70] at (\x,0) {\scriptsize \lab};
    }
    \node[gray!70] at (1,0) {\scriptsize $\times$};
    \node[gray!70] at (2,0) {\scriptsize $\times$};
    \node[gray!70] at (3,0) {\scriptsize $\times$};
    \draw[thick] (1.5,1.7)--(1,0);
    \draw[thick] (1.5,1.7)--(2.5,.85);
    \fill (1.5,1.7) circle (1.5pt);
    \draw[thick] (2.5,.85)--(2,0);
    \draw[thick] (2.5,.85)--(3,0);
    \fill (2.5,.85) circle (1.5pt);
\end{tikzpicture}}

\newcommand{\fTwoCode}{%
\begin{tikzpicture}[x=.46cm,y=.46cm,baseline=-0.6ex,line cap=round,line join=round]
    \draw[gray!65] (0,0)--(3,0);
    \foreach \x/\lab in {0/1,1/2,2/3,3/4} {
        \node[below,gray!70] at (\x,0) {\scriptsize \lab};
    }
    \node[gray!70] at (0,0) {\scriptsize $\times$};
    \node[gray!70] at (1,0) {\scriptsize $\times$};
    \node[gray!70] at (2,0) {\scriptsize $\times$};
    \draw[thick] (1.5,1.7)--(.5,.85);
    \draw[thick] (1.5,1.7)--(2,0);
    \fill (1.5,1.7) circle (1.5pt);
    \draw[thick] (.5,.85)--(0,0);
    \draw[thick] (.5,.85)--(1,0);
    \fill (.5,.85) circle (1.5pt);
\end{tikzpicture}}

\newcommand{\fTop}{%
\begin{tikzpicture}[x=.46cm,y=.46cm,baseline=-0.6ex,line cap=round,line join=round]
    \draw[gray!65] (0,0)--(3,0);
    \foreach \x/\lab in {0/1,1/2,2/3,3/4} {
        \node[below,gray!70] at (\x,0) {\scriptsize \lab};
    }
    \node[gray!70] at (0,0) {\scriptsize $\times$};
    \node[gray!70] at (1,0) {\scriptsize $\times$};
    \node[gray!70] at (2,0) {\scriptsize $\times$};
    \node[gray!70] at (3,0) {\scriptsize $\times$};
    \draw[thick] (1.5,1.7)--(.5,.85);
    \draw[thick] (1.5,1.7)--(2.5,.85);
    \fill (1.5,1.7) circle (1.5pt);
    \draw[thick] (.5,.85)--(0,0);
    \draw[thick] (.5,.85)--(1,0);
    \fill (.5,.85) circle (1.5pt);
    \draw[thick] (2.5,.85)--(2,0);
    \draw[thick] (2.5,.85)--(3,0);
    \fill (2.5,.85) circle (1.5pt);
\end{tikzpicture}}

\begin{tikzpicture}[
    x=1cm,
    y=1cm,
    line cap=round,
    line join=round,
    forestnode/.style={align=center, inner sep=2pt},
    edgelabel/.style={fill=white, rounded corners=1pt, inner sep=1pt, font=\scriptsize}
]

\node[forestnode] (bot) at (0,0)
    {$\emptyforest$};

\node[forestnode] (a) at (-6,3.1)
    {$(0,0,1)$\\[2pt]\fThree};

\node[forestnode] (b) at (0,3.1)
    {$(0,1)$\\[2pt]\fTwo};

\node[forestnode] (c) at (6,3.1)
    {$(1)$\\[2pt]\fOne};

\node[forestnode] (ab) at (-6,6.2)
    {$(0,1,1)$\\[2pt]\fZeroOneOne};

\node[forestnode] (ac) at (0,6.2)
    {$(1,0,1)$\\[2pt]\fOneThree};

\node[forestnode] (bc) at (6,6.2)
    {$(2)$\\[2pt]\fTwoCode};

\node[forestnode] (top) at (0,9.3)
    {$(2,0,1)$\\[2pt]\fTop};

\draw[thick] (bot)-- node[pos=.58,edgelabel,left,xshift=-2pt] {$3$} (a);
\draw[thick] (bot)-- node[pos=.50,edgelabel,left,xshift=-6pt] {$2$} (b);
\draw[thick] (bot)-- node[pos=.58,edgelabel,right,xshift=2pt] {$1$} (c);

\draw[thick] (a)-- node[pos=.55,edgelabel,left,xshift=-2pt] {$2,3$} (ab);
\draw[thick] (a)-- node[pos=.30,edgelabel,below,yshift=-1pt] {$1$} (ac);
\draw[thick] (b)-- node[pos=.34,edgelabel,left,xshift=-2pt] {$3$} (ab);
\draw[thick] (b)-- node[pos=.34,edgelabel,right,xshift=2pt] {$1$} (bc);
\draw[thick] (c)-- node[pos=.30,edgelabel,below,yshift=-1pt] {$3$} (ac);
\draw[thick] (c)-- node[pos=.55,edgelabel,right,xshift=2pt] {$1,2$} (bc);

\draw[thick] (ab)-- node[pos=.24,edgelabel,left,xshift=-2pt] {$1$} (top);
\draw[thick] (ac)-- node[pos=.55,edgelabel,right,xshift=2pt] {$1,2,3$} (top);
\draw[thick] (bc)-- node[pos=.24,edgelabel,right,xshift=2pt] {$3$} (top);

\end{tikzpicture}

\caption{The interval $[\emptyforest,(2,0,1)]$ in the quasisymmetric Bruhat order.}
\label{fig:interval-201}
\endgroup
\end{figure}

A useful first example is the interval whose top has code $(2,0,1)$, drawn in
Figure~\ref{fig:interval-201}. An edge label $i$ records that the upper forest
appears in the expansion of $(x_1+\cdots+x_i)P$ for the lower forest; a label
such as $1,2$ means the cover arises for both $i=1$ and $i=2$.

Each forest in this cube is determined by its support subset. For example,
$(1,0,1)\lessdot(2,0,1)$ fills the missing point between two singleton
components. The classification says that this subset correspondence holds
throughout a lower interval precisely when its top is up-down.

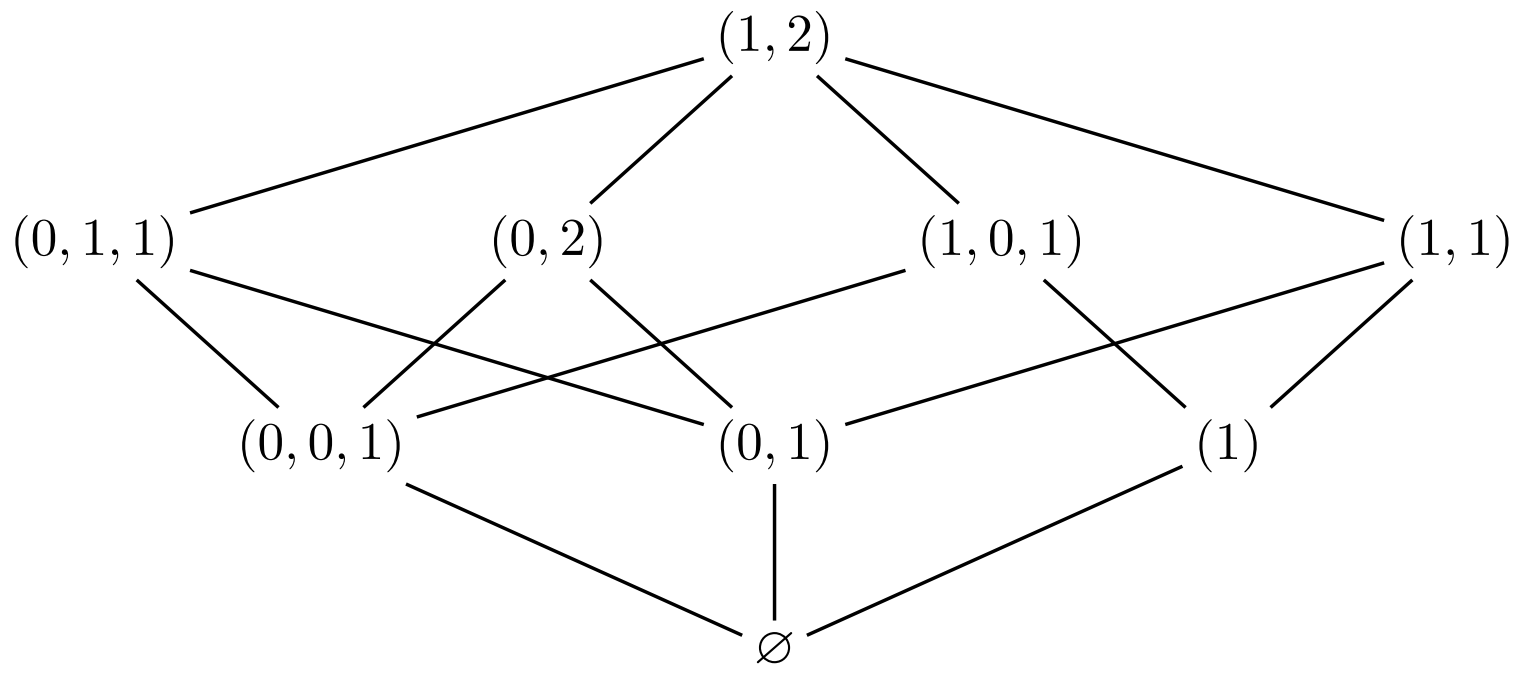
\begin{figure}[H]
\centering
\begin{tikzpicture}[
    x=1cm,
    y=.9cm,
    every node/.style={font=\small},
    vertex/.style={inner sep=2pt},
    edge/.style={line width=.65pt}
]
\node[vertex] (bot) at (0,0) {$\emptyforest$};

\node[vertex] (a) at (-3,1.5) {$(0,0,1)$};
\node[vertex] (b) at (0,1.5) {$(0,1)$};
\node[vertex] (c) at (3,1.5) {$(1)$};

\node[vertex] (d) at (-4.5,3.0) {$(0,1,1)$};
\node[vertex] (e) at (-1.5,3.0) {$(0,2)$};
\node[vertex] (f) at (1.5,3.0) {$(1,0,1)$};
\node[vertex] (g) at (4.5,3.0) {$(1,1)$};

\node[vertex] (top) at (0,4.5) {$(1,2)$};

\draw[edge] (bot)--(a);
\draw[edge] (bot)--(b);
\draw[edge] (bot)--(c);

\draw[edge] (a)--(d);
\draw[edge] (a)--(e);
\draw[edge] (a)--(f);
\draw[edge] (b)--(d);
\draw[edge] (b)--(e);
\draw[edge] (b)--(g);
\draw[edge] (c)--(f);
\draw[edge] (c)--(g);

\draw[edge] (d)--(top);
\draw[edge] (e)--(top);
\draw[edge] (f)--(top);
\draw[edge] (g)--(top);
\end{tikzpicture}
\caption{The interval $[\emptyforest,(1,2)]$, with four elements in rank two.}
\label{fig:interval-12}
\end{figure}

The interval below $(1,2)$ already shows that the poset is not
Eulerian. Figure~\ref{fig:interval-12} draws the
full Hasse diagram. The rank numbers are $(1,3,4,1)$.

\begin{observation}[Non-Eulerian]
\label{obs:non-eulerian}
The quasisymmetric Bruhat order is not Eulerian.
\end{observation}

\begin{proof}
For the interval $[\emptyforest,(1,2)]$, the M\"obius value is
\[
  \mu(\emptyforest,(1,2))
  =-\bigl(1-3+4\bigr)=-2,
\]
which is not $(-1)^3$. Hence the interval is not Eulerian.
\end{proof}

\subsection{The Positive Direction}

\begin{lemma}[One-tree lower covers]
\label{lem:updown-one-tree-covers}
Let $T$ be a single up-down tree. The lower covers of $T$ are 
forests obtained by deleting one support point and rebuilding the remaining
consecutive runs as up-down trees with the inherited turn. In particular,
every lower cover of $T$ is a disjoint union of up-down trees.
\end{lemma}

\begin{proof}
Use the witness-word description of Monk covers. A witness word for an
up-down tree has the form
\[
  \text{an interleaving of }(a_1,\dots,a_p)
  \text{ and }(b_1,\dots,b_q),
  \quad\text{followed by }r.
\]
The final letter $r$ becomes the root. The $a$-sequence builds one straight
branch, and the $b$-sequence builds the other; the two sequences may interleave
before $r$, but the relative order inside each branch is forced.

Figure~\ref{fig:updown-deletion-cases} shows the two representative cases.
For the displayed tree, use the decorated witness word
\[
  w=5^{[1]}\,1^{[1]}\,1^{[2]}\,3^{[1]}\,3^{[2]}.
\]
The final letter $3^{[2]}$ is the root. The upper row deletes this root; the
lower row deletes the branch letter $1^{[1]}$.

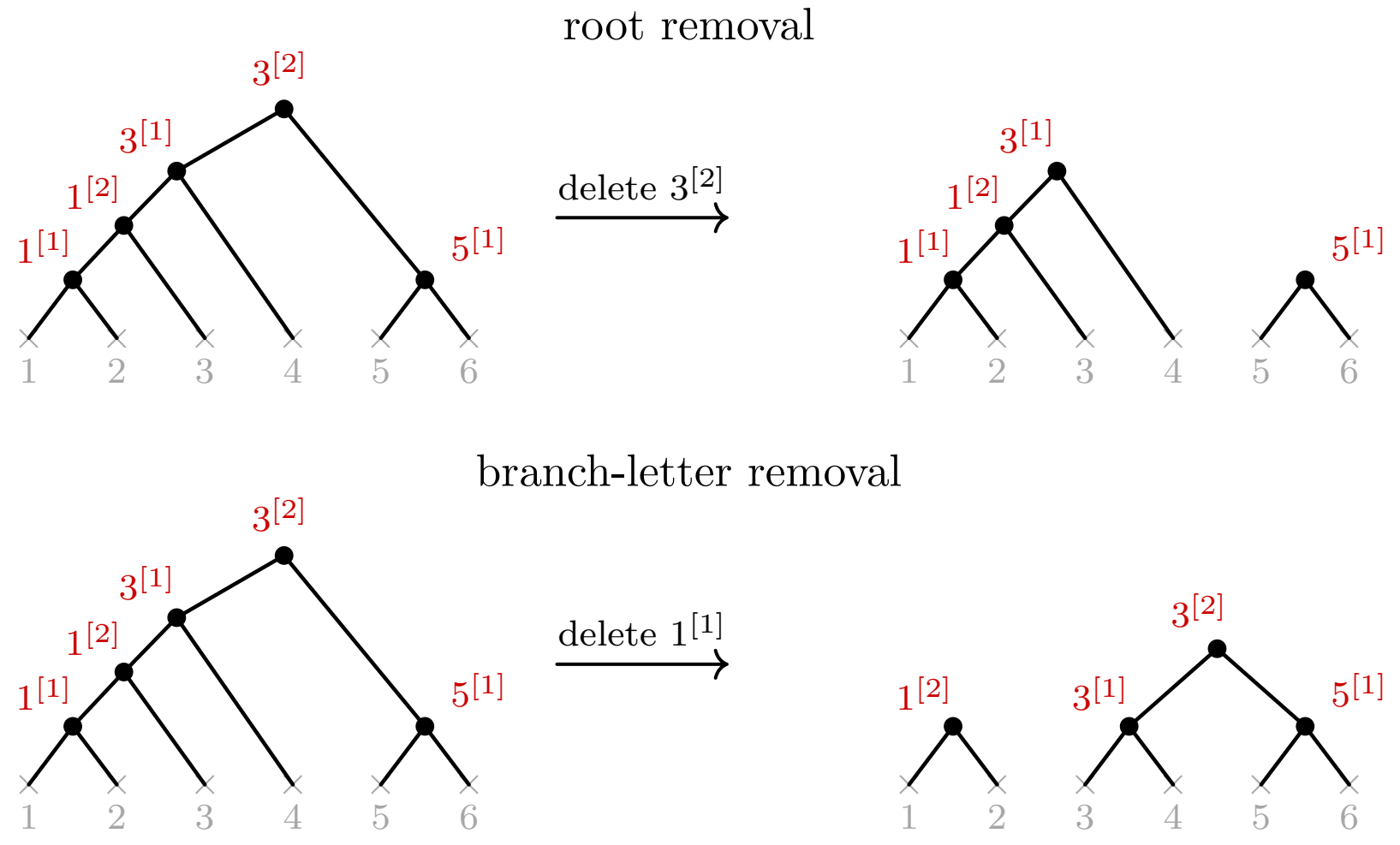
\begin{figure}[H]
\centering
\begin{tikzpicture}[
  x=.68cm,y=.60cm,
  line cap=round,
  line join=round,
  forestedge/.style={line width=.85pt},
  dot/.style={circle,fill=black,inner sep=1.45pt},
  lab/.style={red!80!black,font=\scriptsize,fill=white,inner sep=.7pt},
  leaf/.style={gray!70,font=\scriptsize},
  every node/.style={font=\scriptsize}
]
  \node[font=\small] at (.5,4.05) {root removal};
  \begin{scope}[shift={(-7,0)}]
    \foreach \x/\ell in {0/1,1/2,2/3,3/4,4/5,5/6} {
      \node[leaf] at (\x,0) {$\times$};
      \node[leaf] at (\x,-.42) {$\ell$};
    }
    \coordinate (r1) at (.5,.75);
    \coordinate (r2) at (1.08,1.45);
    \coordinate (r3) at (1.68,2.15);
    \coordinate (r5) at (4.5,.75);
    \coordinate (r4) at (2.9,2.95);
    \draw[forestedge] (r1)--(0,0) (r1)--(1,0);
    \draw[forestedge] (r2)--(r1) (r2)--(2,0);
    \draw[forestedge] (r3)--(r2) (r3)--(3,0);
    \draw[forestedge] (r5)--(4,0) (r5)--(5,0);
    \draw[forestedge] (r4)--(r3) (r4)--(r5);
    \node[dot] at (r1) {}; \node[dot] at (r2) {}; \node[dot] at (r3) {};
    \node[dot] at (r4) {}; \node[dot] at (r5) {};
    \node[lab] at (.18,1.18) {$1^{[1]}$};
    \node[lab] at (.74,1.88) {$1^{[2]}$};
    \node[lab] at (1.35,2.60) {$3^{[1]}$};
    \node[lab] at (2.86,3.48) {$3^{[2]}$};
    \node[lab] at (5.12,1.22) {$5^{[1]}$};
  \end{scope}

  \draw[->,line width=.75pt] (-1.0,1.55)--(.95,1.55)
    node[midway,above,font=\scriptsize] {delete $3^{[2]}$};

  \begin{scope}[shift={(3,0)}]
    \foreach \x/\ell in {0/1,1/2,2/3,3/4,4/5,5/6} {
      \node[leaf] at (\x,0) {$\times$};
      \node[leaf] at (\x,-.42) {$\ell$};
    }
    \coordinate (s1) at (.5,.75);
    \coordinate (s2) at (1.08,1.45);
    \coordinate (s3) at (1.68,2.15);
    \coordinate (s5) at (4.5,.75);
    \draw[forestedge] (s1)--(0,0) (s1)--(1,0);
    \draw[forestedge] (s2)--(s1) (s2)--(2,0);
    \draw[forestedge] (s3)--(s2) (s3)--(3,0);
    \draw[forestedge] (s5)--(4,0) (s5)--(5,0);
    \node[dot] at (s1) {}; \node[dot] at (s2) {}; \node[dot] at (s3) {};
    \node[dot] at (s5) {};
    \node[lab] at (.18,1.18) {$1^{[1]}$};
    \node[lab] at (.74,1.88) {$1^{[2]}$};
    \node[lab] at (1.35,2.60) {$3^{[1]}$};
    \node[lab] at (5.12,1.22) {$5^{[1]}$};
  \end{scope}

  \begin{scope}[shift={(0,-5.75)}]
    \node[font=\small] at (.5,4.05) {branch-letter removal};
    \begin{scope}[shift={(-7,0)}]
      \foreach \x/\ell in {0/1,1/2,2/3,3/4,4/5,5/6} {
        \node[leaf] at (\x,0) {$\times$};
        \node[leaf] at (\x,-.42) {$\ell$};
      }
      \coordinate (m1) at (.5,.75);
      \coordinate (m2) at (1.08,1.45);
      \coordinate (m3) at (1.68,2.15);
      \coordinate (m5) at (4.5,.75);
      \coordinate (m4) at (2.9,2.95);
      \draw[forestedge] (m1)--(0,0) (m1)--(1,0);
      \draw[forestedge] (m2)--(m1) (m2)--(2,0);
      \draw[forestedge] (m3)--(m2) (m3)--(3,0);
      \draw[forestedge] (m5)--(4,0) (m5)--(5,0);
      \draw[forestedge] (m4)--(m3) (m4)--(m5);
      \node[dot] at (m1) {}; \node[dot] at (m2) {}; \node[dot] at (m3) {};
      \node[dot] at (m4) {}; \node[dot] at (m5) {};
      \node[lab] at (.18,1.18) {$1^{[1]}$};
      \node[lab] at (.74,1.88) {$1^{[2]}$};
      \node[lab] at (1.35,2.60) {$3^{[1]}$};
      \node[lab] at (2.86,3.48) {$3^{[2]}$};
      \node[lab] at (5.12,1.22) {$5^{[1]}$};
    \end{scope}

    \draw[->,line width=.75pt] (-1.0,1.55)--(.95,1.55)
      node[midway,above,font=\scriptsize] {delete $1^{[1]}$};

    \begin{scope}[shift={(3,0)}]
      \foreach \x/\ell in {0/1,1/2,2/3,3/4,4/5,5/6} {
        \node[leaf] at (\x,0) {$\times$};
        \node[leaf] at (\x,-.42) {$\ell$};
      }
      \coordinate (t1) at (.5,.75);
      \coordinate (t3) at (2.5,.75);
      \coordinate (t5) at (4.5,.75);
      \coordinate (t4) at (3.5,1.75);
      \draw[forestedge] (t1)--(0,0) (t1)--(1,0);
      \draw[forestedge] (t3)--(2,0) (t3)--(3,0);
      \draw[forestedge] (t5)--(4,0) (t5)--(5,0);
      \draw[forestedge] (t4)--(t3) (t4)--(t5);
      \node[dot] at (t1) {}; \node[dot] at (t3) {};
      \node[dot] at (t4) {}; \node[dot] at (t5) {};
      \node[lab] at (.18,1.18) {$1^{[2]}$};
      \node[lab] at (2.18,1.18) {$3^{[1]}$};
      \node[lab] at (3.30,2.25) {$3^{[2]}$};
      \node[lab] at (5.12,1.22) {$5^{[1]}$};
    \end{scope}
  \end{scope}
\end{tikzpicture}
\caption{Root and branch-letter deletions in an up-down tree.}
\label{fig:updown-deletion-cases}
\end{figure}

If the deleted letter is $r$, the remaining word is just an interleaving of the
two branch words, so the two branches are no longer joined. The result is a
disjoint union of a pure up tree and a pure down tree.

Now suppose the deleted letter is $a_k$ on one branch. The letters before
$a_k$ still build the initial straight part of that branch. After the deletion,
each later branch letter either attaches to this initial part or is the first
letter forced to face the missing support point. That first nonattaching letter
becomes terminal, and from it upward the remaining letters have the same
up-down form as before. Hence the deletion splits the branch into a straight
initial tree and a later up-down tree. The same argument applies to the
$b$-branch by reflection.

Thus every lower cover has the stated form. Conversely,
Theorem~\ref{thm:support-removing} gives a lower cover removing any prescribed
support point $x$. The preceding argument identifies it with the rebuilt
forest, so every forest in the statement occurs.
\end{proof}

For a fixed up-down forest $G$ and a subset $A\subseteq\supp(G)$, define
$\UD_G(A)$ as follows. On each component support of $G$, keep exactly the
points of $A$ in that component support. Each remaining consecutive run is
rebuilt as an up-down tree using the turn inherited from the original
component of $G$.

\begin{lemma}[Lower covers of up-down forests]
\label{lem:updown-forest-covers}
If $G$ is an up-down forest, then its lower covers are exactly
\[
  \UD_G(\supp(G)\setminus\{x\}),\qquad x\in\supp(G).
\]
\end{lemma}

\begin{proof}
By Lemma~\ref{lem:covers-are-component-covers}, a lower cover changes only the
component containing the removed support point $x$.
Lemma~\ref{lem:updown-one-tree-covers} identifies the result as
$\UD_G(\supp(G)\setminus\{x\})$. Conversely, its one-tree cover remains a
cover after adjoining the unchanged components, by
Lemma~\ref{lem:component-cover-remain}.
\end{proof}

\begin{theorem}[Up-down forests give Boolean intervals]
\label{thm:updown-boolean-forward}
If $G$ is an up-down forest and $S=\supp(G)$, then
\[
  [\emptyforest,G]\cong B_{|S|},
\]
the Boolean lattice on $S$.
\end{theorem}

\begin{proof}
By Lemma~\ref{lem:updown-forest-covers}, descending from $G$ successively
deletes support points and rebuilds the remaining runs. Hence every element
below $G$ is $\UD_G(A)$ for some $A\subseteq S$, and every such forest occurs
by deleting $S\setminus A$. The same lemma gives exactly the covers
\[
  \UD_G(A\setminus\{x\})\lessdot\UD_G(A),\qquad x\in A.
\]
Thus $\UD_G(A)\mapsto A$ is a bijection carrying covers to single-element
inclusions, as required.
\end{proof}

\subsection{Zigzag Isolation}

A path of internal vertices determines a zigzag shape by replacing every
off-path branch with a leaf. We show that this shape occurs as a forest below
the original one: choose its witness letters so that erasing the surrounding
letters preserves the path.

\begin{theorem}[Zigzag isolation]
\label{thm:zigzag-isolation}
Let $F$ be an indexed forest, and let $Z$ be the zigzag shape of a downward
path of internal vertices in $F$. Then a positively indexed copy of $Z$ lies
below $F$ in the quasisymmetric Bruhat order.
\end{theorem}

\begin{proof}
List the chosen vertices from bottom to top, and let $i$ be the inorder
support label of the bottom vertex. First build the off-path subtrees below
the top vertex using their canonical postorder words, as in the proof of
Theorem~\ref{thm:support-removing}. Insert the bottom vertex with value $i$.
For each successive parent, use value $i$ if the existing path is its left
child, and value $i-1$ if it is its right child. Choose fresh decorations
increasing along each run of right growth and decreasing along each run of
left growth.

These choices rebuild the prescribed path with all its attached subtrees.
Indeed, its current root and the next off-path piece are adjacent in the
root list. A piece on the left has root value at most $i-2$, or is a bare
leaf at most $\underline{i-1}$; a piece on the right has root value at least
$i+1$, or is a bare leaf at least $\underline{i+1}$. The chosen parent lies
between the two entries and joins precisely these pieces. The initial
bottom insertion joins its two children by the same root-list update.
Once this rooted subtree is built, complete the remaining postorder
insertions to obtain a witness word $w$ for $F$.

Now retain only the chosen path letters, obtaining a subword $z$. Its first
letter creates a singleton at $i$. Each later letter is on the same side of
the current root as before, so it attaches a new parent on the prescribed
side, with a bare leaf replacing the omitted subtree. Thus $\OmegaF(z)$
has shape $Z$. More explicitly, if there are $L$ left-growth steps and $R$
right-growth steps, then
\[
  \supp(\OmegaF(z))=[i-L,i+R].
\]
The $L$ left-growth ancestors have distinct positive inorder labels below
$i$, so this support is positive. Finally, restoring the omitted letters at
their original positions gives $\OmegaF(z)\leq\OmegaF(w)=F$ by
Corollary~\ref{cor:restoration}.
\end{proof}

For the Boolean classification, only a three-vertex path with a change of
direction is needed. Isolating it gives the small non-Boolean interval used
below.

\subsection{Necessity}

We now prove the converse: if the lower interval below a forest is Boolean,
then the forest must be up-down.

\begin{proposition}[Boolean forests have Boolean components]
\label{prop:boolean-components}
If $F=F_1\sqcup\cdots\sqcup F_r$ and $[\emptyforest,F]$ is Boolean, then each
lower interval $[\emptyforest,F_j]$ is Boolean.
\end{proposition}

\begin{proof}
This follows immediately from the component product decomposition
(Proposition~\ref{prop:component-product}). Each $[\emptyforest,F_j]$ is a
factor, hence also an interval inside the Boolean lattice $[\emptyforest,F]$.
Intervals inside Boolean lattices are Boolean.
\end{proof}

\begin{definition}[Bent corner]
\label{def:bent-corner}
A \emph{bent corner} in a tree is a rooted path of three internal vertices
\[
  p\longrightarrow r\longrightarrow c
\]
where the direction changes at the middle vertex $r$. Thus either $r$ is the
right child of $p$ and $c$ is the left child of $r$, or $r$ is the left child
of $p$ and $c$ is the right child of $r$.
\end{definition}

\begin{lemma}
\label{lem:nonupdown-has-bend}
A connected tree is up-down if and only if it has no bent corner.
\end{lemma}

\begin{proof}
In an up-down tree, every downward path of internal vertices keeps the same
left/right direction, so no bent corner occurs.

Conversely, without a bent corner, an internal left child cannot have an
internal right child, and an internal right child cannot have an internal
left child. Hence the root's left subtree is a straight left branch and its
right subtree is a straight right branch. This is precisely an up-down tree.
\end{proof}

\begin{lemma}[Bent corners force the rank-three obstruction]
\label{lem:bend-forces-obstruction}
If a connected tree $T$ is not up-down, then there is a forest $H\leq T$ such
that $H$ is a shifted copy of $(1,2)$ or $(2,1)$.
\end{lemma}

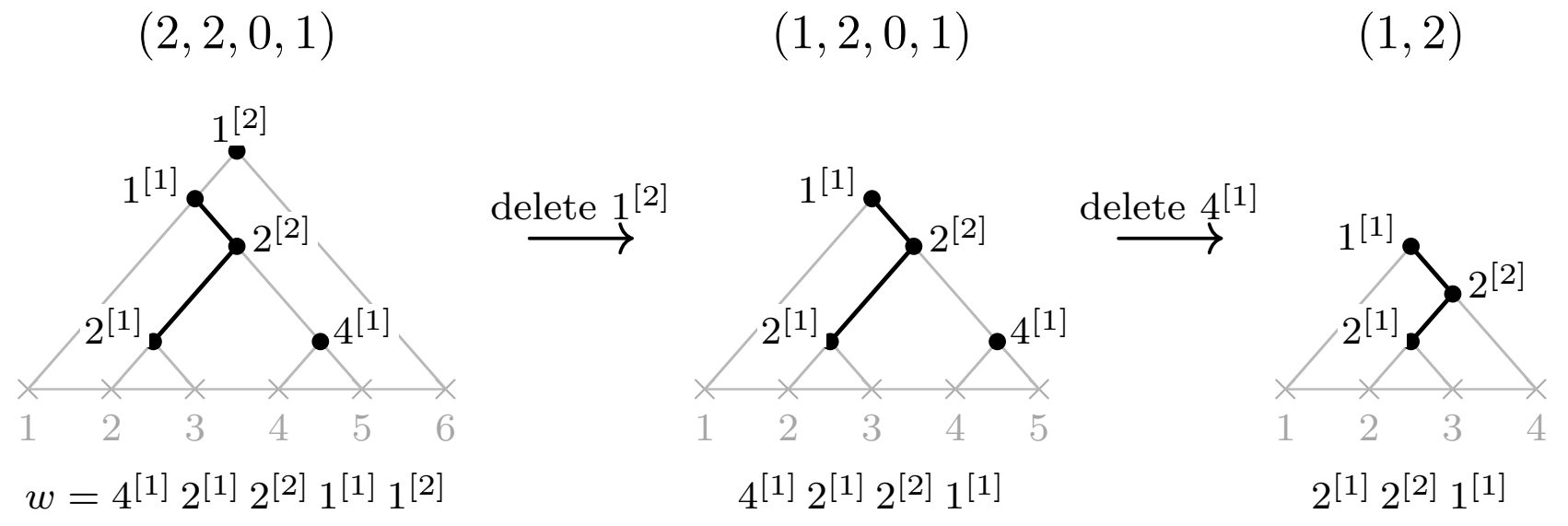
\begin{figure}[H]
\centering
\begin{tikzpicture}[
    x=.60cm,
    y=.57cm,
    line cap=round,
    line join=round,
    focusedge/.style={black,line width=.9pt},
    lightedge/.style={gray!55,line width=.55pt},
    numberline/.style={gray!55,line width=.45pt},
    dot/.style={circle,fill=black,inner sep=1.25pt},
    lab/.style={font=\scriptsize,fill=white,inner sep=.7pt},
    leaf/.style={gray!70,font=\scriptsize},
    word/.style={font=\scriptsize,align=center},
    title/.style={font=\small}
]
\begin{scope}[shift={(0,0)}]
    \node[title] at (2.5,4.45) {$(2,2,0,1)$};
    \draw[numberline] (0,0)--(5,0);
    \foreach \x/\ell in {0/1,1/2,2/3,3/4,4/5,5/6} {
        \node[leaf] at (\x,0) {$\times$};
        \node[leaf] at (\x,-.48) {$\ell$};
    }
    \coordinate (sRoot) at (2.5,3.0);
    \coordinate (sP) at (2.0,2.4);
    \coordinate (sR) at (2.5,1.8);
    \coordinate (sC) at (1.5,.6);
    \coordinate (sExtra) at (3.5,.6);
    \draw[lightedge] (sRoot)--(sP) (sRoot)--(5,0);
    \draw[lightedge] (sP)--(0,0);
    \draw[lightedge] (sC)--(1,0) (sC)--(2,0);
    \draw[lightedge] (sExtra)--(3,0) (sExtra)--(4,0);
    \draw[lightedge] (sR)--(sExtra);
    \draw[focusedge] (sP)--(sR) (sR)--(sC);
    \foreach \v in {sRoot,sP,sR,sC,sExtra} {\node[dot] at (\v) {};}
    \node[lab] at (2.55,3.34) {$1^{[2]}$};
    \node[lab] at (1.48,2.58) {$1^{[1]}$};
    \node[lab] at (3.05,1.96) {$2^{[2]}$};
    \node[lab] at (1.04,.80) {$2^{[1]}$};
    \node[lab] at (4.02,.80) {$4^{[1]}$};
    \node[word] at (2.5,-1.28)
        {$w=4^{[1]}\,2^{[1]}\,2^{[2]}\,1^{[1]}\,1^{[2]}$};
\end{scope}

\draw[->,line width=.7pt] (6.0,1.9)--(7.25,1.9)
    node[midway,above,font=\scriptsize] {delete $1^{[2]}$};

\begin{scope}[shift={(8.1,0)}]
    \node[title] at (2,4.45) {$(1,2,0,1)$};
    \draw[numberline] (0,0)--(4,0);
    \foreach \x/\ell in {0/1,1/2,2/3,3/4,4/5} {
        \node[leaf] at (\x,0) {$\times$};
        \node[leaf] at (\x,-.48) {$\ell$};
    }
    \coordinate (mP) at (2.0,2.4);
    \coordinate (mR) at (2.5,1.8);
    \coordinate (mC) at (1.5,.6);
    \coordinate (mExtra) at (3.5,.6);
    \draw[lightedge] (mP)--(0,0);
    \draw[lightedge] (mC)--(1,0) (mC)--(2,0);
    \draw[lightedge] (mExtra)--(3,0) (mExtra)--(4,0);
    \draw[lightedge] (mR)--(mExtra);
    \draw[focusedge] (mP)--(mR) (mR)--(mC);
    \foreach \v in {mP,mR,mC,mExtra} {\node[dot] at (\v) {};}
    \node[lab] at (1.48,2.58) {$1^{[1]}$};
    \node[lab] at (3.05,1.96) {$2^{[2]}$};
    \node[lab] at (1.04,.80) {$2^{[1]}$};
    \node[lab] at (4.02,.80) {$4^{[1]}$};
    \node[word] at (2,-1.28)
        {$4^{[1]}\,2^{[1]}\,2^{[2]}\,1^{[1]}$};
\end{scope}

\draw[->,line width=.7pt] (13.05,1.9)--(14.3,1.9)
    node[midway,above,font=\scriptsize] {delete $4^{[1]}$};

\begin{scope}[shift={(15.05,0)}]
    \node[title] at (1.5,4.45) {$(1,2)$};
    \draw[numberline] (0,0)--(3,0);
    \foreach \x/\ell in {0/1,1/2,2/3,3/4} {
        \node[leaf] at (\x,0) {$\times$};
        \node[leaf] at (\x,-.48) {$\ell$};
    }
    \coordinate (tP) at (1.5,1.8);
    \coordinate (tR) at (2.0,1.2);
    \coordinate (tC) at (1.5,.6);
    \draw[lightedge] (tP)--(0,0);
    \draw[lightedge] (tR)--(3,0);
    \draw[lightedge] (tC)--(1,0) (tC)--(2,0);
    \draw[focusedge] (tP)--(tR) (tR)--(tC);
    \foreach \v in {tP,tR,tC} {\node[dot] at (\v) {};}
    \node[lab] at (.98,1.98) {$1^{[1]}$};
    \node[lab] at (2.55,1.38) {$2^{[2]}$};
    \node[lab] at (1.04,.80) {$2^{[1]}$};
    \node[word] at (1.5,-1.28)
        {$2^{[1]}\,2^{[2]}\,1^{[1]}$};
\end{scope}
\end{tikzpicture}
\caption{Isolating a bent corner from $(2,2,0,1)$}
\label{fig:bent-corner-isolation}
\end{figure}

\begin{proof}
By Lemma~\ref{lem:nonupdown-has-bend}, the tree contains a bent corner
\[
  p\longrightarrow r\longrightarrow c.
\]
These three vertices form a path-shaped zigzag pattern rooted at $p$ inside
$T$. By zigzag isolation (Theorem~\ref{thm:zigzag-isolation}), the three-vertex
zigzag with this turn word lies below $T$ in the quasisymmetric Bruhat order as a standalone
tree $H$. Depending on whether the turn is right-then-left or left-then-right,
the code of $H$ is a shifted copy of $(1,2)$ or $(2,1)$.
Figure~\ref{fig:bent-corner-isolation} shows this isolation for the tree
$(2,2,0,1)$, whose bent corner is pared down to the shifted copy of $(1,2)$ one
Monk cover at a time.
\end{proof}

\begin{lemma}[The rank-three obstruction is not Boolean]
\label{lem:rank-three-obstruction}
The interval $[\emptyforest,(1,2)]$ is not Boolean. The same is true for
$[\emptyforest,(2,1)]$ and for every shifted copy of either interval.
\end{lemma}

\begin{proof}
As drawn in Figure~\ref{fig:interval-12}, the interval below $(1,2)$ has rank numbers
\[
  1,\ 3,\ 4,\ 1.
\]
A Boolean lattice of rank $3$ has rank numbers $1,3,3,1$, so this interval is
not Boolean. The interval below $(2,1)$ is the reflection of this one, and
shifted copies have isomorphic lower intervals by
Proposition~\ref{prop:symmetries}.
\end{proof}

\begin{proof}[Proof of Theorem~\ref{thm:boolean-classification}]
The forward construction is Theorem~\ref{thm:updown-boolean-forward}.
Conversely, if $[\emptyforest,F]$ is Boolean, then every component interval
is Boolean by Proposition~\ref{prop:boolean-components}. A component that is
not up-down would contain a shifted rank-three obstruction $H$ below it, by
Lemma~\ref{lem:bend-forces-obstruction}. Then $[\emptyforest,H]$ would be an
interval of a Boolean lattice, contradicting
Lemma~\ref{lem:rank-three-obstruction}. Thus every component is up-down.
\end{proof}

\begin{remark}
\label{rem:lattice-classification}
The same obstruction distinguishes lattices. In
Figure~\ref{fig:interval-12} the atoms $(0,0,1)$ and $(0,1)$ have two
minimal common upper bounds, $(0,1,1)$ and $(0,2)$, so
$[\emptyforest,(1,2)]$ is not a lattice. Every interval of a
lattice is a lattice, and Lemma~\ref{lem:bend-forces-obstruction} embeds
one of these intervals below every forest that is not up-down. Combined
with Theorem~\ref{thm:boolean-classification}, the lower interval
$[\emptyforest,F]$ is therefore a lattice if and only if $F$ is an
up-down forest.
\end{remark}

\section{Shellability}
\label{sec:shellability}

Boolean lower intervals are shellable. The next question is whether this
property persists for the non-Boolean intervals. We prove it for the hook
family $(n,1)$ by ordering coatoms recursively, then use code-raising to
reduce the general conjecture to gap-free codes.

\subsection{Order Complexes and Shellability}

For a finite poset $P$, the \emph{order complex} $\Delta(P)$ is the simplicial
complex whose faces are the chains of $P$. When $P$ has a minimum $\hat0$ and a
maximum $\hat1$, the complex of interest is that of the \emph{proper part}
$P\setminus\{\hat0,\hat1\}$, so shellability of a bounded interval $[x,y]$ means
shellability of the complex of chains strictly between $x$ and $y$. Every lower
interval $[\emptyforest,F]$ is bounded, with $\hat0=\emptyforest$ and
$\hat1=F$, and it is graded by Section~\ref{sec:first-structure}; its maximal
chains therefore all have length $|\supp(F)|$ and its order complex is pure.

A \emph{shelling} of a pure simplicial complex is an ordering of its facets
$C_1,C_2,\dots$ such that, for each $j>1$, the intersection of $C_j$ with the
union of the preceding facets is pure of codimension one in $C_j$. In an
order complex, the facets are maximal chains. A shelling therefore assembles
the maximal chains along their codimension-one faces. The usual empty-face
convention applies in dimension zero.

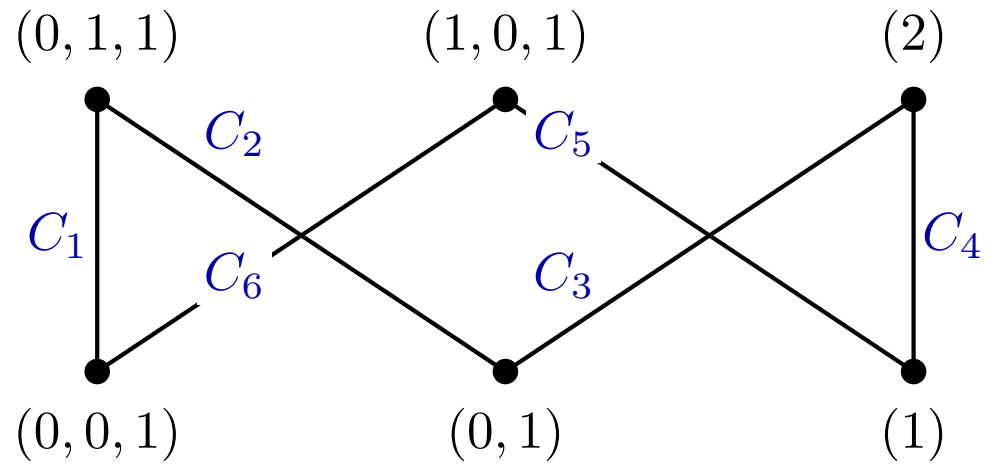
\begin{figure}[H]
\centering
\begin{tikzpicture}[
  x=1cm,y=1cm,
  point/.style={circle,fill=black,inner sep=1.7pt},
  cover/.style={black,line width=.8pt},
  chain/.style={font=\small,text=blue!65!black,fill=white,inner sep=1.2pt},
  code/.style={font=\small}
]
  \coordinate (a) at (-2.7,0);
  \coordinate (b) at (0,0);
  \coordinate (c) at (2.7,0);
  \coordinate (ab) at (-2.7,1.8);
  \coordinate (ac) at (0,1.8);
  \coordinate (bc) at (2.7,1.8);
  \draw[cover] (a)--node[chain,left] {$C_1$} (ab);
  \draw[cover] (b)--node[chain,pos=.76,above right] {$C_2$} (ab);
  \draw[cover] (b)--node[chain,pos=.24,above left] {$C_3$} (bc);
  \draw[cover] (c)--node[chain,right] {$C_4$} (bc);
  \draw[cover] (c)--node[chain,pos=.76,above left] {$C_5$} (ac);
  \draw[cover] (a)--node[chain,pos=.24,above right] {$C_6$} (ac);
  \foreach \p in {a,b,c,ab,ac,bc} \node[point] at (\p) {};
  \node[code,below=3pt] at (a) {$(0,0,1)$};
  \node[code,below=3pt] at (b) {$(0,1)$};
  \node[code,below=3pt] at (c) {$(1)$};
  \node[code,above=3pt] at (ab) {$(0,1,1)$};
  \node[code,above=3pt] at (ac) {$(1,0,1)$};
  \node[code,above=3pt] at (bc) {$(2)$};
\end{tikzpicture}
\caption{Shelling $[\emptyforest,(2,0,1)]$ in the order $C_1,\ldots,C_6$.}
\label{fig:interval-201-shelling}
\end{figure}

Figure~\ref{fig:interval-201-shelling} uses the Boolean interval from
Figure~\ref{fig:interval-201}. Its proper part is a six-cycle; the order
complex is built by adding its six edges in order.

\begin{fact}[{\cite[Corollary~3.1.9(b)]{Wac07}}]
\label{fact:interval-heredity}
If a poset $P$ is shellable, then every interval $[x,y]\subseteq P$ is
shellable.
\end{fact}

Fact~\ref{fact:interval-heredity} allows us to deduce shellability of
subintervals from shellability of a larger interval.

\begin{remark}
\label{rem:boolean-shellable}
If $G$ is an up-down forest, then $[\emptyforest,G]$ is the Boolean lattice on
$\supp(G)$ by Theorem~\ref{thm:boolean-classification}, so the standard
labeling of $B_{|\supp(G)|}$ by the element added at each cover is an
EL-labeling~\cite[Section~3.2.1]{Wac07}; Figure~\ref{fig:interval-201-shelling}
illustrates a shelling in rank $3$. When
$|\supp(G)|\ge 2$, the proper part is the barycentric subdivision of a simplex
boundary, so the open interval is a sphere of dimension
$|\supp(G)|-2$~\cite[Example~1.1.1]{Wac07}, and
$\mu(\emptyforest,G)=(-1)^{|\supp(G)|}$~\cite[Exercise~1.2.3]{Wac07}. The M\"obius value $-2$ recorded for
$(1,2)$ in Observation~\ref{obs:non-eulerian} already rules this out below the
smallest trees that are not up-down.
\end{remark}

\subsection{Recursive Coatom Orderings and the Hook Family}

A coatom is an element immediately below the top. Ordering coatoms divides
the maximal chains into blocks according to their last step; within each
block, the same construction is applied to a smaller interval. The following
conditions ensure that the blocks attach in a shelling order.

\begin{definition}[Recursive coatom ordering]
\label{def:rco}
Let $P$ be a bounded graded poset with minimum $\hat0$ and maximum $\hat1$. An
ordering $q_1,\dots,q_t$ of the coatoms of $P$ is a \emph{recursive coatom
ordering} if $\rk(P)\le 2$, or if both of the following hold:
\begin{enumerate}[label=(R\arabic*),leftmargin=2.8em]
  \item for each $j$, the interval $[\hat0,q_j]$ has a recursive coatom
  ordering in which the coatoms lying below some earlier $q_i$ ($i<j$) come
  first;
  \item if $i<j$ and $y\le q_i,q_j$, then there are $\ell<j$ and
  $z\lessdot q_j$ with $y\le z\le q_\ell$.
\end{enumerate}
The conditions are automatic in rank $2$, so this is equivalent to the dual of
Wachs's length-$1$ definition~\cite[Definition~4.2.1]{Wac07}; rank $0$ uses the empty order.
A recursive coatom ordering is equivalent to dual CL-shellability, and hence
implies shellability~\cite[Theorem~3.2]{BW83}.
\end{definition}

We apply this to the hook intervals. Write $S_m$ for the straight tree with
code $(m)$ and $T_m$ for the hook tree with code $(m,1)$, and set
\[
  \mathcal{B}_m=[\emptyforest,S_m],
  \qquad
  \mathcal{I}_m=[\emptyforest,T_m].
\]
The straight tree $S_m$ is up-down, so $\mathcal{B}_m$ is the Boolean lattice of
rank $m$ by Theorem~\ref{thm:updown-boolean-forward}. The hook $T_m$ is not
up-down for $m\ge 2$, so the $\mathcal{I}_m$ are the first non-Boolean intervals
we shell. Throughout we use the component product decomposition of
Proposition~\ref{prop:component-product}, which identifies
$[\emptyforest,F_1\sqcup F_2]$ with
$[\emptyforest,F_1]\times[\emptyforest,F_2]$ when the two components have
separated supports. Translates with leftmost support position $a$ are written
$S_m[a]$ and $T_m[a]$, and empty components are omitted.

\begin{lemma}[Hook coatoms]
\label{lem:hook-coatoms}
The coatoms of $\mathcal{I}_n$ are
\begin{align*}
  D_{k,n}&=S_k[1]\sqcup S_{n-k}[k+2], &&1\le k\le n,\\
  C_{0,n}&=S_n[2],\\
  C_{k,n}&=T_k[1]\sqcup S_{n-k-1}[k+3], &&1\le k\le n-1.
\end{align*}
In code coordinates,
\[
  D_{k,n}=(k,0^k,n-k),
  \qquad
  C_{0,n}=(0,n),
  \qquad
  C_{k,n}=(k,1,0^k,n-k-1).
\]
\end{lemma}

\begin{figure}[H]
\centering
\begingroup
\newcommand{\HookLeafLine}{%
  \draw[gray!65,line width=.45pt] (0,0)--(3,0);
  \foreach \x/\lab in {0/1,1/2,2/3,3/4} {
    \node[gray!65] at (\x,0) {\scriptsize $\times$};
    \node[below=2pt,font=\scriptsize,text=gray!65] at (\x,0) {\lab};
  }
}
\newcommand{\HookTop}{%
  \begin{tikzpicture}[x=.43cm,y=.43cm,baseline=-.6ex]
    \HookLeafLine
    \draw[thick] (2,2.3)--(.8,1.5)--(0,0);
    \draw[thick] (2,2.3)--(3,0);
    \draw[thick] (.8,1.5)--(1.5,.7)--(1,0);
    \draw[thick] (1.5,.7)--(2,0);
    \fill (2,2.3) circle (1.4pt);
    \fill (.8,1.5) circle (1.4pt);
    \fill (1.5,.7) circle (1.4pt);
  \end{tikzpicture}}
\newcommand{\HookCZero}{%
  \begin{tikzpicture}[x=.43cm,y=.43cm,baseline=-.6ex]
    \HookLeafLine
    \draw[thick] (2.5,1.5)--(1.5,.7)--(1,0);
    \draw[thick] (1.5,.7)--(2,0);
    \draw[thick] (2.5,1.5)--(3,0);
    \fill (2.5,1.5) circle (1.4pt); \fill (1.5,.7) circle (1.4pt);
  \end{tikzpicture}}
\newcommand{\HookDOne}{%
  \begin{tikzpicture}[x=.43cm,y=.43cm,baseline=-.6ex]
    \HookLeafLine
    \draw[thick] (0,0)--(.5,.7)--(1,0);
    \draw[thick] (2,0)--(2.5,.7)--(3,0);
    \fill (.5,.7) circle (1.4pt); \fill (2.5,.7) circle (1.4pt);
  \end{tikzpicture}}
\newcommand{\HookCOne}{%
  \begin{tikzpicture}[x=.43cm,y=.43cm,baseline=-.6ex]
    \HookLeafLine
    \draw[thick] (.5,1.5)--(0,0);
    \draw[thick] (.5,1.5)--(1.5,.7)--(1,0);
    \draw[thick] (1.5,.7)--(2,0);
    \fill (.5,1.5) circle (1.4pt); \fill (1.5,.7) circle (1.4pt);
  \end{tikzpicture}}
\newcommand{\HookDTwo}{%
  \begin{tikzpicture}[x=.43cm,y=.43cm,baseline=-.6ex]
    \HookLeafLine
    \draw[thick] (1.5,1.5)--(.5,.7)--(0,0);
    \draw[thick] (.5,.7)--(1,0);
    \draw[thick] (1.5,1.5)--(2,0);
    \fill (1.5,1.5) circle (1.4pt); \fill (.5,.7) circle (1.4pt);
  \end{tikzpicture}}
\begin{tikzpicture}[
  x=1cm,y=1cm,
  forestnode/.style={align=center,inner sep=3pt,font=\small}
]
  \node[forestnode] (top) at (0,2.5) {$T_2=(2,1)$\\[3pt]\HookTop};
  \node[forestnode] (C0) at (-3.9,0) {$C_{0,2}=(0,2)$\\[3pt]\HookCZero};
  \node[forestnode] (D1) at (-1.3,0) {$D_{1,2}=(1,0,1)$\\[3pt]\HookDOne};
  \node[forestnode] (C1) at (1.3,0) {$C_{1,2}=(1,1)$\\[3pt]\HookCOne};
  \node[forestnode] (D2) at (3.9,0) {$D_{2,2}=(2)$\\[3pt]\HookDTwo};
  \draw[line width=.8pt] (C0)--(top) (D1)--(top) (C1)--(top) (D2)--(top);
\end{tikzpicture}
\endgroup
\caption{The coatoms of the hook $T_2=(2,1)$.}
\label{fig:hook-coatom-example}
\end{figure}
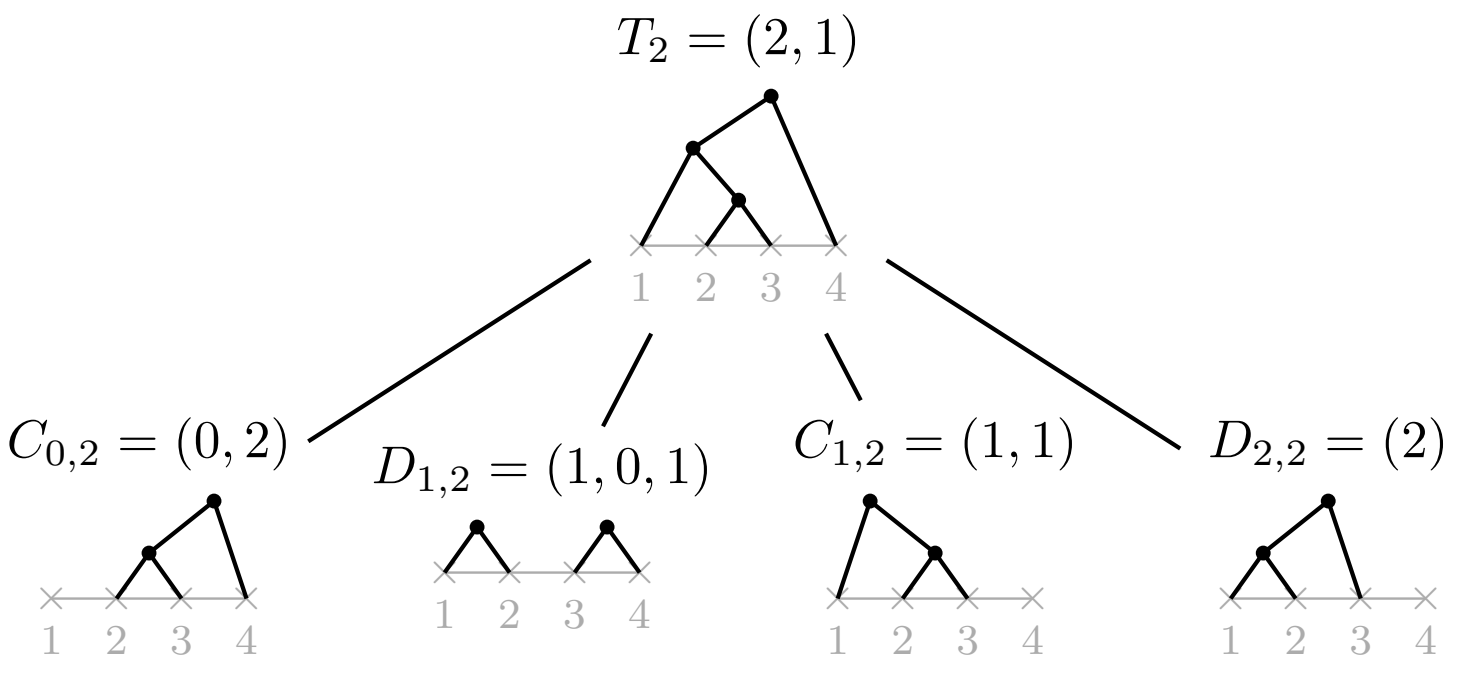

\begin{proof}
Deleting a letter either straightens the hook or leaves a smaller hook
beside a straight component. To distinguish these outcomes, track the
single missing support point as the remaining letters are inserted.

Let $w=a_1\cdots a_{n+1}$ witness $T_n$. Its internal path forces
bottom-to-top insertion, since insertion preserves existing children. Thus
\[
 \underline2<a_1<\underline3,\qquad
 \underline1<a_2<a_1,\qquad
 a_{t-1}<a_t<\underline{t+1}\quad(3\le t\le n+1).
\]
By Lemma~\ref{lem:word-monk}, every lower cover comes from deleting some
$a_j$ from such a witness. By Theorem~\ref{thm:subword-support} and the
grading, its support is $[1,n+1]\setminus\{p\}$ for some $p$.

If $j=1$, the remaining word is increasing. Each insertion starts a
component or extends the rightmost component to the right; hence its
components are straight trees on the consecutive runs of the support.
The result is $C_{0,n}$ for $p=1$, and $D_{p-1,n}$ otherwise.

Suppose $j\ge3$. The surviving prefix is $T_{j-2}$ on $[1,j-1]$, with
gap $j$. Until a separate right component appears, just before a later
letter $a_t$ the lower prefix is a hook on $[1,t-2]$, whereas the upper
prefix has support $[1,t-1]$. The lower root, initially $a_{j-1}$ and
subsequently $a_{t-1}$, is smaller than $a_t$, and its next bare entry
is $\underline t$. If $a_t<\underline t$, insertion extends the hook,
fills $t-1$, and moves the gap to $t$. Otherwise it starts a singleton
at $t$, leaving gap $t-1$. In this second case the upper component and
the rightmost lower component now share root $a_t$ and next bare entry
$\underline{t+2}$. Every later letter lies between their common root
and right frontier, so it extends the rightmost lower component and
the gap stops moving. A final gap $p$ therefore gives
\[
 T_{p-2}[1]\sqcup S_{n+1-p}[p+1]=C_{p-2,n},
 \qquad 3\le p\le n+1.
\]
An empty right component includes deletion of the final letter and the
case when the gap never stops moving.

For $j=2$, the case $n=1$ gives $C_{0,1}$. Otherwise, if $a_3>a_1$,
the shortened word is increasing and begins at position $2$, so
position $1$ remains absent and the result is $C_{0,n}$. If $a_3<a_1$,
then $\underline1<a_2<a_3<a_1$: inserting $a_3$ joins $\underline1$ to
$a_1$, recreating $T_1$ on $[1,2]$. The gap is now $3$, and the same
moving-gap argument applies.

Conversely, the witness $2,1,3,\ldots,n+1$ gives $C_{0,n}$ by deleting
$1$, and $C_{k,n}$ by deleting $k+2$ for $1\le k<n$. Here integers
denote decorated letters, with decoration $1$ for distinct values and
increasing decorations on repetitions. For $1\le k\le n$, the witness
\[
 2,\underbrace{1,\ldots,1}_{k\text{ terms}},
   \underbrace{k+2,\ldots,k+2}_{n-k\text{ terms}}
\]
also builds $T_n$; deleting $2$ gives $D_{k,n}$. The same support and
word-cover results make these deletions covers. Finally, distinct
omitted positions distinguish the listed forests except possibly
$C_{k,n}$ and $D_{k+1,n}$, whose left components $T_k$ and $S_{k+1}$
differ. Reading off their codes gives the stated coordinates.
\end{proof}

Order the coatoms of $\mathcal{I}_n$ by
\[
  \mathcal{O}_n:\quad
  C_{0,n},\,D_{1,n},\,C_{1,n},\,D_{2,n},\,\dots,\,
  D_{n-1,n},\,C_{n-1,n},\,D_{n,n}.
\]
By Lemma~\ref{lem:hook-coatoms} and Proposition~\ref{prop:component-product},
the intervals below the coatoms factor as
\begin{align*}
  [\emptyforest,D_{k,n}]&\cong\mathcal{B}_k\times\mathcal{B}_{n-k}
      &&(1\le k\le n),\\
  [\emptyforest,C_{0,n}]&\cong\mathcal{B}_n,\\
  [\emptyforest,C_{k,n}]&\cong\mathcal{I}_k\times\mathcal{B}_{n-k-1}
      &&(1\le k<n).
\end{align*}
For a coatom $q$, let $U_{<q}$ denote the union of the order ideals generated by
the coatoms preceding $q$ in $\mathcal{O}_n$; we call an element \emph{old} at
$q$ if it lies in $U_{<q}$.

\begin{lemma}[Old intersections]
\label{lem:old-intersections}
For $1\le k\le n$,
\[
  [\emptyforest,D_{k,n}]\cap U_{<D_{k,n}}
  =(\mathcal{B}_k\setminus\{S_k\})\times\mathcal{B}_{n-k},
\]
and for $1\le k<n$,
\[
  [\emptyforest,C_{k,n}]\cap U_{<C_{k,n}}
  =(\mathcal{I}_k\setminus\{T_k\})\times\mathcal{B}_{n-k-1}.
\]
\end{lemma}

\begin{proof}
In each product, an element is old exactly when its first factor is
below the top. Take $(A,B)\le D_{k,n}$ under the identification
$[\emptyforest,D_{k,n}]\cong\mathcal{B}_k\times\mathcal{B}_{n-k}$. If $A<S_k$,
choose a position $p\in[1,k]$ missing from $A$; enlarging $A$ to the Boolean
coatom that omits only $p$ places $(A,B)$ below $C_{0,n}$ when $p=1$ and below
$D_{p-1,n}$ when $p>1$, and either coatom precedes $D_{k,n}$, so $(A,B)$ is old.
If $A=S_k$, support containment excludes every earlier coatom except possibly
$C_{k-1,n}$ (when $k\geq2$), whose first component $T_{k-1}$ shares the support of $S_k$ but is
a distinct forest of the same rank, hence incomparable to $S_k$. This proves
the first identity.

Now take $(A,B)\le C_{k,n}$ under
$[\emptyforest,C_{k,n}]\cong\mathcal{I}_k\times\mathcal{B}_{n-k-1}$. If
$A<T_k$, choose a coatom $q$ of $\mathcal{I}_k$ above $A$; if $q=C_{a,k}$ then
$(A,B)$ lies below $C_{a,n}$, and if $q=D_{a,k}$ then it lies below $D_{a,n}$,
both of which precede $C_{k,n}$. Conversely, every coatom preceding $C_{k,n}$
omits a support position in $[1,k+1]$, so no element with first coordinate
$T_k$ lies below it. This proves the second identity.
\end{proof}

\begin{lemma}[Products of recursive coatom orderings]
\label{lem:rco-product}
If $p_1,\ldots,p_s$ and $q_1,\ldots,q_t$ are recursive coatom orderings
of finite bounded graded posets $P$ and $Q$, then
\[
 (p_1,\hat1_Q),\ldots,(p_s,\hat1_Q),
 (\hat1_P,q_1),\ldots,(\hat1_P,q_t)
\]
is a recursive coatom ordering of $P\times Q$.
\end{lemma}

\begin{proof}
Induct on the sum of the ranks; total rank at most two is immediate, and
a rank-zero factor reduces to the other factor. Write
$A_i=(p_i,\hat1_Q)$ and $B_j=(\hat1_P,q_j)$.

Below $A_i$, the old coatoms are exactly $(r,\hat1_Q)$ with
$r\lessdot p_i$ lying below an earlier $p_h$. Choose the (R1) ordering
below $p_i$ and apply induction with that factor first.

Below $B_j$, all $(p_i,q_j)$ are old. The remaining old coatoms are
$(\hat1_P,r)$ with $r\lessdot q_j$ lying below an earlier $q_h$.
The (R1) ordering below $q_j$, followed by product induction with $P$
first, therefore puts exactly the old coatoms first. This proves (R1).

For (R2), pairs of the same type use the corresponding factor's (R2).
For a mixed pair, every common lower bound $x\leq A_i,B_j$ satisfies
\[
  x\leq(p_i,q_j)\lessdot B_j,
  \qquad (p_i,q_j)\leq A_i,
\]
and $A_i$ is earlier. Thus (R2) holds.
\end{proof}

\begin{theorem}
\label{thm:n1-shellable}
For every $n\ge 1$, the order $\mathcal{O}_n$ is a recursive coatom ordering of
$\mathcal{I}_n$. Therefore $[\emptyforest,(n,1)]$ is dual CL-shellable and
shellable.
\end{theorem}

\begin{proof}
Induct on $n$. The interval $\mathcal{I}_1$ is Boolean of rank $2$, so every
ordering of its coatoms is recursive. Assume $\mathcal{O}_k$ is a recursive
coatom ordering of $\mathcal{I}_k$ for every $k<n$.

For (R1), $C_{0,n}$ has a Boolean lower interval and no earlier coatom,
so any coatom order there suffices. At every later coatom, the old
coatoms are exactly those from the first factor, by
Lemma~\ref{lem:old-intersections}. Put this block first. The resulting
order is recursive: for a $D$-type coatom both factors are Boolean, and
for $C_{k,n}$ use $\mathcal O_k$ on the first factor and
Lemma~\ref{lem:rco-product}.

For (R2), suppose $y\le q$ also lies below an earlier coatom of
$\mathcal{O}_n$. Lemma~\ref{lem:old-intersections} places $y$ below a
first-factor coatom $z\lessdot q$, and the proof of that lemma exhibits $z$
below an earlier member of $\mathcal{O}_n$. This verifies (R2) and completes the
induction.
\end{proof}
\subsection{Reduction to Zigzag Tops}

Every forest lies below a zigzag forest, obtained by filling the gaps in
its code. Shellability of the zigzag lower intervals would therefore settle
the general case.

\begin{definition}[Zigzag forest]
\label{def:zigzag-forest}
A forest $Z$ is a \emph{zigzag forest} if, after deleting trailing zeros from
its code, every remaining entry is positive:
\[
  \code(Z)=(a_1,\dots,a_k),\qquad a_j>0\ \text{for all }j.
\]
\end{definition}

These are the zigzag forests of Nadeau--Spink--Tewari~\cite[Section~3.5]{NST26},
normalized to have $1$ in their support, together with $\emptyforest$.
Their nonempty members are precisely the trees whose internal vertices form a
single ancestor chain and whose support begins at $1$~\cite[Proposition~3.10]{NT23}.
Translation preserves principal lower intervals up to isomorphism by
Proposition~\ref{prop:symmetries}. The hook $T_2=(2,1)$ in
Figure~\ref{fig:hook-coatom-example} is one example.

\begin{proposition}[Zigzag completion]
\label{prop:zigzag-completion}
For every forest $F$ there is a zigzag forest $Z(F)$ with $F\leq Z(F)$, obtained
by filling the zero positions of $\code(F)$ one cover at a time.
\end{proposition}

\begin{proof}
Replace each zero before the last nonzero code entry by $1$, one at a
time. Each replacement is a Monk cover by
Lemma~\ref{lem:code-raising}. The final code has no gaps, giving the
required $Z(F)$ above $F$. For the empty forest, take
$Z(\emptyforest)=\emptyforest$.
\end{proof}

\begin{corollary}[Reduction to zigzag forests]
\label{cor:zigzag-reduction}
If $[\emptyforest,Z]$ is shellable for every zigzag forest $Z$, then
$[\emptyforest,F]$ is shellable for every forest $F$.
\end{corollary}

\begin{proof}
Given $F$, Proposition~\ref{prop:zigzag-completion} gives
$\emptyforest\leq F\leq Z(F)$, so $[\emptyforest,F]$ is an interval inside
$[\emptyforest,Z(F)]$. If the latter is shellable, then
Fact~\ref{fact:interval-heredity} makes the former shellable.
\end{proof}

\subsection{Computational Evidence and the Conjecture}

A recursive-coatom search found orderings for the lower intervals below
every binary tree with $m$ internal vertices and leaves $1,\dots,m+1$, for
$m=4,5,6$. There are respectively $14$, $42$, and $132$ such tops. In all
$188$ cases, the induced maximal-chain ordering passed a separate shelling
check. These finite computations support the conjecture but do not establish
it in general.

By Corollary~\ref{cor:zigzag-reduction}, shellability of the zigzag lower
intervals would settle every principal lower interval.

\begin{conjecture}[Shellability]
\label{conj:zigzag-shellable}
\label{conj:lower-shellable}
Every zigzag lower interval $[\emptyforest,Z]$ is shellable; consequently every
lower interval $[\emptyforest,F]$ is shellable.
\end{conjecture}

\section*{Acknowledgments}

I am grateful to my mentor Dora Woodruff for her guidance, feedback, and
mentorship throughout this project. I also thank the PRIMES-USA program for
making this research possible.

\end{document}